\documentclass[11pt]{article}
\usepackage{amssymb,amsfonts,amsmath,curves}
\usepackage{epsfig}
\usepackage{graphicx}
\usepackage{color}
\usepackage{anysize}
\usepackage{hyperref}

\marginsize{2cm}{2cm}{1cm}{1cm}

\begin{document}
\baselineskip=14pt

\newtheorem{defin}{Definition}[section]
\newtheorem{Prop}{Proposition}
\newtheorem{teo}{Theorem}[section]
\newtheorem{ml}{Main Lemma}
\newtheorem{con}{Conjecture}
\newtheorem{cond}{Condition}
\newtheorem{conj}{Conjecture}
\newtheorem{prop}[teo]{Proposition}
\newtheorem{lem}{Lemma}[section]
\newtheorem{rmk}[teo]{Remark}
\newtheorem{cor}{Corollary}[section]

\newcommand{\beq}{\begin{equation}}
\newcommand{\eeq}{\end{equation}}
\newcommand{\beqn}{\begin{eqnarray}}
\newcommand{\beqnn}{\begin{eqnarray*}}
\newcommand{\eeqn}{\end{eqnarray}}
\newcommand{\eeqnn}{\end{eqnarray*}}
\newcommand{\bprop}{\begin{prop}}
\newcommand{\eprop}{\end{prop}}
\newcommand{\bt}{\begin{teo}}
\newcommand{\bcor}{\begin{cor}}
\newcommand{\ecor}{\end{cor}}
\newcommand{\bcon}{\begin{con}}
\newcommand{\econ}{\end{con}}
\newcommand{\bcond}{\begin{cond}}
\newcommand{\econd}{\end{cond}}
\newcommand{\bconj}{\begin{conj}}
\newcommand{\econj}{\end{conj}}
\newcommand{\et}{\end{teo}}
\newcommand{\brm}{\begin{rmk}}
\newcommand{\erm}{\end{rmk}}
\newcommand{\blem}{\begin{lem}}
\newcommand{\elem}{\end{lem}}
\newcommand{\ben}{\begin{enumerate}}
\newcommand{\een}{\end{enumerate}}
\newcommand{\bei}{\begin{itemize}}
\newcommand{\eei}{\end{itemize}}
\newcommand{\bdf}{\begin{defin}}
\newcommand{\edf}{\end{defin}}
\newcommand{\bpr}{\begin{proof}}
\newcommand{\epr}{\end{proof}}

\newenvironment{proof}{\noindent {\em Proof}.\,\,}{\hspace*{\fill}$\halmos$\medskip}

\newcommand{\halmos}{\rule{1ex}{1.4ex}}
\def \qed {{\hspace*{\fill}$\halmos$\medskip}}

\newcommand{\fr}{\frac}
\newcommand{\Z}{{\mathbb Z}}
\newcommand{\R}{{\mathbb R}}
\newcommand{\E}{{\mathbb E}}
\newcommand{\C}{{\mathbb C}}
\renewcommand{\P}{{\mathbb P}}
\newcommand{\N}{{\mathbb N}}
\newcommand{\Q}{{\mathbb Q}}
\newcommand{\var}{{\mathbb V}}
\renewcommand{\S}{{\cal S}}
\newcommand{\T}{{\cal T}}
\newcommand{\W}{{\cal W}}
\newcommand{\X}{{\cal X}}
\newcommand{\Y}{{\cal Y}}
\newcommand{\h}{{\cal H}}
\newcommand{\Fi}{{\cal F}}
\newcommand{\Ni}{{\cal N}}

\renewcommand{\a}{\alpha}
\renewcommand{\b}{\beta}
\newcommand{\g}{\gamma}
\newcommand{\G}{\Gamma}
\renewcommand{\L}{\Lambda}
\renewcommand{\l}{\lambda}
\renewcommand{\d}{\delta}
\newcommand{\D}{\Delta}
\newcommand{\e}{\epsilon}
\newcommand{\eps}{\epsilon}
\newcommand{\s}{\sigma}
\newcommand{\B}{{\cal B}}
\renewcommand{\o}{\omega}

\newcommand{\nn}{\nonumber}
\renewcommand{\=}{&=&}
\renewcommand{\>}{&>&}
\newcommand{\<}{&<&}
\renewcommand{\le}{\leq}
\newcommand{\+}{&+&}

\newcommand{\pa}{\partial}
\newcommand{\ffrac}[2]{{\textstyle\frac{{#1}}{{#2}}}}
\newcommand{\dif}[1]{\ffrac{\partial}{\partial{#1}}}
\newcommand{\diff}[1]{\ffrac{\partial^2}{{\partial{#1}}^2}}
\newcommand{\difif}[2]{\ffrac{\partial^2}{\partial{#1}\partial{#2}}}

\title{Survival Probability of a Random Walk Among a Poisson System of Moving Traps}
\author{Alexander Drewitz$^{\,1}$ \and J\"urgen G\"artner$^{\,2}$ \and Alejandro F.\ Ram\'irez$^{\,3}$ \and Rongfeng Sun$^{\,4}$}
\date{Nov 1, 2011}

\maketitle

\footnotetext[1]{Departement Mathematik, Eidgen\"ossische Technische Hochschule Z\"urich, R\"amistrasse 101, 8092 Z\"urich, Switzerland. Email: \mbox{alexander.drewitz@math.ethz.ch}}

\footnotetext[2]{Institut f\"ur Mathematik, Technische Universit\"at Berlin, Sekr.~MA 7-5, Str.~des 17.~Juni 136, 10623 Berlin, Germany. Email: jg@mail.math.tu-berlin.de}

\footnotetext[3]{Facultad de Matem\'aticas, Pontificia Universidad Cat\'olica de Chile, Vicu$\rm \tilde{n}$a Mackenna 4860, Macul, Santiago, Chile. Email: aramirez@mat.puc.cl}

\footnotetext[4]{Department of Mathematics, National University of Singapore, 10 Lower Kent Ridge Road, 119076 Singapore. Email:
matsr@nus.edu.sg}

\begin{abstract}

We review some old and prove some new results on the survival probability of a random
walk among a Poisson system of moving traps on $\Z^d$, which can also be
interpreted as the solution of a parabolic Anderson model with a random
time-dependent potential. We show that the annealed survival probability
decays asymptotically as $e^{-\lambda_1\sqrt{t}}$ for $d=1$, as $e^{-\lambda_2 t/\log t}$
for $d=2$, and as $e^{-\lambda_dt}$ for $d\geq 3$, where $\lambda_1$ and $\lambda_2$ can be identified
explicitly. In addition, we show that the quenched survival probability decays asymptotically
as $e^{-\tilde \lambda_d t}$, with $\tilde \lambda_d>0$ for all $d\geq 1$. A key ingredient in bounding
the annealed survival probability is what is known in the physics literature as the Pascal
principle, which asserts that the annealed survival probability is maximized
if the random walk stays at a fixed position. A corollary of independent
interest is that the expected cardinality of the range of a continuous time
symmetric random walk increases under perturbation by a deterministic path.


\medskip

\end{abstract}


\section{Introduction}

\subsection{Model and results}
Let $X:=(X(t))_{t\geq 0}$ be a simple symmetric random walk on $\Z^d$ with jump rate $\kappa \geq 0$, and let $(Y^y_j)_{1\leq j\leq N_y, y\in\Z^d}$ be a collection
of independent simple symmetric random walks on $\Z^d$ with jump rate $\rho> 0$, where $N_y$ is the number of walks that start at each $y\in\Z^d$ at
time $0$, $(N_y)_{y\in\Z^d}$ are i.i.d.\ Poisson distributed with mean $\nu>0$, and $Y^y_j:=(Y^y_j(t))_{t\geq 0}$ denotes the $j$-th walk starting at $y$ at time $0$.
Let us denote the number of walks $Y$ at position $x\in\Z^d$ at time $t\geq 0$ by
\beq\label{xi}
\xi(t,x) := \sum_{y\in \Z^d, 1\leq j\leq N_y} \delta_x(Y^y_j(t)).
\eeq
It is easy to see that for each $t\geq 0$, $(\xi(t, x))_{x\in\Z^d}$ are i.i.d.\ Poisson distributed with mean $\nu$, so that $(\xi(t,\cdot))_{t\geq 0}$ is a stationary
process, and furthermore it is reversible in the sense that $(\xi(t, \cdot))_{0\leq t\leq T}$ is equally distributed with $(\xi(T-t, \cdot))_{0\leq t\leq T}$.
We will interpret the collection of walks $Y$ as traps, and at each time $t$, the walk $X$ is killed with rate $\gamma \xi(t, X(t))$ for some parameter
$\gamma>0$. Conditional on the realization of the field of traps $\xi$, the probability that the walk $X$ survives by time $t$ is given by
\beq\label{quenchedZ}
Z^{\gamma}_{t,\xi} := \E^X_0\Big[\exp\Big\{-\gamma \int_0^t \xi(s,X(s)) \, {\rm d}s\Big\}\Big],
\eeq
where $\E^X_0$ denotes expectation with respect to $X$ with $X(0)=0$. We call this the {\em quenched survival probability}, which depends on the random medium $\xi$.
When we furthermore average over $\xi$, which we denote by $\E^\xi$, we obtain the {\em annealed survival probability}
\beq\label{annealedZ}
\E^\xi[Z^\gamma_{t,\xi}] = \E^\xi\E^X_0\Big[\exp\Big\{-\gamma \int_0^t \xi(s,X(s))\, {\rm d}s\Big\}\Big].
\eeq
We will study the long time behavior of the annealed and quenched survival probabilities, and in particular, identify their rate of decay and their dependence
on the spatial dimension $d$ and the parameters $\kappa, \rho, \nu$ and $\gamma$.
\bigskip

Here are our main results on the decay rate of the annealed and quenched survival probabilities.
\bt\label{T:annealed}{\bf [Annealed survival probability]}
Assume that $\gamma\in (0,\infty]$, $\kappa\geq 0$, $\rho>0$ and $\nu>0$, then
\beq\label{annealed}
\E^\xi[Z^\gamma_{t,\xi}]\ =\ \left\{
\begin{aligned}
\exp\Big\{-\nu \sqrt{\frac{8\rho t}{\pi}}(1+o(1))\Big\}, & & \qquad d=1, \\
\exp\Big\{-\nu\pi\rho \frac{t}{\log t}(1+o(1))\Big\}, & & \qquad d=2,  \\
\exp\Big\{-\lambda_{d,\gamma,\kappa,\rho, \nu}\, t(1+o(1))\Big\}, & & \qquad d\geq 3,
\end{aligned}
\right.
\eeq
where $\lambda_{d,\gamma,\kappa,\rho, \nu}$ depends on $d$, $\gamma$, $\kappa$, $\rho$, $\nu$, and is called the annealed Lyapunov exponent. Furthermore,
$\lambda_{d,\gamma,\kappa,\rho, \nu}\geq \lambda_{d,\gamma,0,\rho, \nu}=\nu\gamma/(1+\frac{\gamma G_d(0)}{\rho})$, where $G_d(0):=\int_0^\infty p_t(0)\, {\rm d}t$ is the Green function of a simple symmetric random walk on $\Z^d$ with jump rate $1$ and transition kernel $p_t(\cdot)$.
\et
Note that in dimensions $1$ and $2$, the annealed survival probability decays sub-exponentially, and the pre-factor in front of the decay rate is surprisingly independent of $\gamma\in (0,\infty]$ and $\kappa\geq 0$. The key ingredient in the proof is what is known in the physics literature as the {\em Pascal principle}, which asserts that in (\ref{annealedZ}), if we condition on the random walk trajectory $X$, then the annealed survival probability is maximized when $X\equiv 0$. The discrete time version of the Pascal principle was proved by Moreau, Oshanin, B\'enichou and Coppey in~\cite{MOBC03, MOBC04}. We will include the proof for the reader's convenience.
As a corollary of the Pascal principle, we will show in Corollary \ref{C:range} that the expected cardinality of the range of a continuous time symmetric random walk increases under perturbation by a deterministic path.
\medskip

In contrast to the annealed case, the quenched survival probability always decays exponentially.
\bt\label{T:quenched}{\bf [Quenched survival probability]}
Assume that $d\geq 1$, $\gamma>0$, $\kappa\geq 0$, $\rho>0$ and $\nu>0$. Then there exists deterministic $\tilde\lambda_{d,\gamma,\kappa,\rho, \nu}$ depending on
$d, \gamma,\kappa,\rho, \nu$, called the quenched Lyapunov exponent, such that $\P^\xi$-a.s.,
\beq\label{quenched}
Z^\gamma_{t,\xi} = \exp\big\{-\tilde\lambda_{d,\gamma,\kappa,\rho, \nu}\,t(1+o(1))\big\} \quad  \mbox{as } t\to\infty.
\eeq
Furthermore, $0<\tilde\lambda_{d,\gamma,\kappa,\rho, \nu}\leq \gamma\nu+\kappa$ for all $d\geq 1, \gamma>0, \kappa\geq 0, \rho>0$ and $\nu>0$.
\et
{\bf Remark.} When $\gamma<0$, $Z^\gamma_{t,\xi}$ can be interpreted as the expected number of branching random walks in the catalytic medium $\xi$. See Section~\ref{S:lit} for more discussion on this model. As will be outlined at the end of Section~\ref{S:shape}, (\ref{quenched}) also holds
in this case, and lies in the interval $[-\gamma\nu-\kappa, \infty)$.
\medskip

In Proposition~\ref{P:panneal} below, we will also give an upper bound of the same order as in Theorem~\ref{T:annealed} for the survival probability $\E^\xi[Z^\gamma_{t,\xi}]$, where $(\xi(0,x))_{x\in\Z^d}$ is deterministic and satisfies some constraints. These constraints hold asymptotically
a.s.\ for i.i.d.\ Poisson distributed $(\xi(0,x))_{x\in\Z^d}$. Therefore we call this a {\em semi-annealed} bound, which we will use in Section~\ref{S:quenched1} to obtain
sub-exponential bounds on the quenched survival probability in dimensions $1$ and $2$.

\subsection{Relation to the parabolic Anderson model}
The annealed and quenched survival probabilities $Z^\gamma_{t,\xi}$ and $\E^\xi[Z^\gamma_{t,\xi}]$ are closely related to the solution of the parabolic Anderson model (PAM), namely, the solution of the following parabolic equation with random potential $\xi$:
\beq \label{PAM}
\begin{aligned}
\frac{\partial}{\partial t} u(t, x) &= \kappa \Delta u(t, x) - \gamma\, \xi(t, x)\, u(t, x),  \\
u(0,x) &= 1,
\end{aligned}
\qquad \qquad x\in \Z^d,\ t\geq 0,
\eeq
where $\gamma, \kappa$ and $\xi$ are as before, and $\Delta f(x) = \frac{1}{2d} \sum_{\Vert y-x\Vert =1} (f(y)-f(x))$ is the
discrete Laplacian on $\Z^d$, which is also the generator of a simple symmetric random walk on $\Z^d$ with jump rate $1$.

By the Feynman-Kac formula, the solution $u$ admits the representation
\beq\label{FKrep}
u(t,0) = \E^X_0\left[\exp\left\{-\gamma\int_0^t \xi(t-s, X(s)) \, {\rm d}s\right\} \right],
\eeq
which differs from $Z^\gamma_{t,\xi}$ in (\ref{quenchedZ}) by a time reversal in $\xi$. When we average $u(t,0)$ over the random field $\xi$,
by the reversibility of $(\xi(t,\cdot))_{0\leq s\leq t}$, we have
\beq\label{annealedequiv}
\E^\xi[u(t,0)] = \E^\xi\E^X_0\!\!\left[\exp\!\left\{\!-\gamma\!\int_0^t\! \xi(t-s, X(s)) \, {\rm d}s\!\right\}\!\right]\! =\! \E^X_0\E^\xi\!\!\left[\exp\!\left\{\!-\gamma\!\int_0^t \!\xi(s, X(s))\, {\rm d}s\!\right\}\! \right] = \E^\xi[Z^\gamma_{t,\xi}].
\eeq
Therefore Theorem~\ref{T:annealed} also applies  to the annealed solution $\E^\xi[u(t,0)]$. Despite the difference between $Z^\gamma_{t,\xi}$ and $u(t,0)$ due to
time reversal, Theorem~\ref{T:quenched} also holds with $u(t,0)$ in place of $Z^\gamma_{t,\xi}$.

\bt\label{T:quenchedPAM}{\bf [Quenched solution of PAM]} Let $d\geq 1$, $\gamma>0$, $\kappa\geq 0$, $\rho>0$, $\nu>0$ and $\tilde\lambda_{d,\gamma,\kappa,\rho, \nu}>0$ be the same as in Theorem~\ref{T:quenched}. Then $\P^\xi$-a.s.,
\beq\label{quenchedPAM}
u(t,0) = \exp\big\{- \tilde\lambda_{d,\gamma,\kappa,\rho, \nu}\, t(1+o(1))\big\}  \quad \mbox{as } t\to\infty.
\eeq
\et
{\bf Remark.} By Theorem~\ref{T:quenched} and the remark following it, for any $\gamma\in\R$, $t^{-1}\log u(t,0)$ converges in probability to $- \tilde\lambda_{d,\gamma,\kappa,\rho, \nu}$ because $u(t,0)$ is equally distributed with $Z^\gamma_{t,\xi}$. However, we were only able to strengthen this  to almost sure convergence for the $\gamma>0$ case, but not for $\gamma<0$. For a broader investigation of the case $\gamma<0$, see G\"artner, den Hollander, and Maillard~\cite{GHM11}, which is also contained in the present volume.

\subsection{Review of related results}\label{S:lit}
The study of trapping problems has a long history in the mathematics and physics literature. We review some models and results that are most relevant to our problem.

\subsubsection{Immobile Traps} Extensive studies have been carried out for the case of immobile traps, i.e., $\rho=0$ and $\xi(t,\cdot)\equiv \xi(0,\cdot)$ for all $t\geq 0$.
A continuum version is Brownian motion among Poissonian obstacles, where a ball of size $1$ is placed and centered at each point of a mean density $1$ homogeneous Poisson point process in $\R^d$, acting as traps or obstacles, and an independent Brownian motion starts at the origin and is killed at rate $\gamma$ times the number of obstacles it is contained in. Using a large deviation principle for the Brownian motion occupation time measure, Donsker and Varadhan~\cite{DV75} showed that the annealed survival probability decays asymptotically as $\exp\{-C_{d,\gamma}t^{\frac{d}{d+2}}(1+o(1))\}$. Using spectral techniques, Sznitman~\cite{S98} later developed a coarse graining method, known as the {\em method of enlargement of obstacles}, to show that the quenched survival probability decays asymptotically as $\exp\{-\bar C_{d,\gamma}\frac{t}{(\log t)^{2/d}}(1+o(1))\}$. Similar results have also been obtained for random walks among immobile Bernoulli traps (i.e.\ $\xi(0,x)\in\{0,1\}$), see e.g.\ \cite{DV79, B94, A94, A95}. Traps with a more general form of the trapping potential $\xi$ have also been studied in the context of the parabolic Anderson model (see e.g.\ Biskup and K\"onig~\cite{BK01}), where alternative techniques to the method of enlargement of obstacles were developed and the order of sub-exponential decay of the survival probabilities may vary depending on the distribution of $\xi$. Compared to our results in Theorems~\ref{T:annealed} and \ref{T:quenched}, we note that when the traps are moving, both the annealed and quenched survival probabilities decay faster than when the traps are immobile. The heuristic reason is that, the walk survives by finding large space-time regions void of traps, which are easily destroyed if the traps are moving. Another example is a Brownian motion among Poissonian obstacles where the obstacles move with a deterministic drift. It has been shown that the annealed and quenched survival probabilities decay exponentially if the drift is sufficiently large, see e.g.~\cite[Thms.~5.4.7 and 5.4.9]{S98}.

\subsubsection{Mobile Traps} The model we consider here has in fact been studied earlier by Redig in \cite{R94}, where he considered a trapping potential $\xi$ generated by a reversible Markov process, such as a Poisson system of random walks, or the symmetric exclusion process in equilibrium. Using spectral techniques applied to the
process of moving traps viewed from the random walk, he established an exponentially decaying upper bound for the annealed survival probability
when the empirical distribution of the trapping potential, $\frac{1}{t}\int_0^t \xi(s,0)\, {\rm d}s$, satisfies a large deviation principle with scale $t$. This applies
for instance to $\xi$ generated from either a Poisson system of independent random walks or the symmetric exclusion process in equilibrium, in dimensions $d\geq 3$.

\subsubsection{Annihilating Two-type Random Walks} In \cite{BL91}, Bramson and Lebowitz studied a model from chemical physics, where there are two types of particles, $A$s
and $B$s, both starting initially with an i.i.d.\ Poisson distribution on $\Z^d$ with density $\rho_A(0)$ resp.\ $\rho_B(0)$. All particles perform independent simple symmetric random walk with jump rate $1$, particles of the same type do not interact, and when two particles of opposite types meet, they annihilate each other. This system models a chemical reaction $A+B\to {\rm inert}$. It was shown in \cite{BL91} that when $\rho_A(0)=\rho_B(0)>0$, then $\rho_A(t)$ and $\rho_B(t)$ (the densities of the $A$ and $B$ particles at time $t$) decay with the order $t^{-d/4}$ in dimensions $1\leq d\leq 4$, and decay with the order $t^{-1}$ in $d\geq 4$. When $\rho_A(0)>\rho_B(0)>0$, it was shown that $\rho_A(t)\to \rho_A(0)-\rho_B(0)$ as $t\to\infty$, and $-\log\rho_B(t)$ increases with the order $\sqrt{t}$ in $d=1$, $t/\log t$ in $d=2$, and $t$ in $d\geq 3$, which is the same as in Theorem~\ref{T:annealed}.  Heuristically, as $\rho_B(t)\to 0$ and $\rho_A(t)\to \rho_A(0)-\rho_B(0)>0$, we can effectively model the $B$ particles as uncorrelated single random walks among a Poisson field of moving traps with density $\rho_A(0)-\rho_B(0)$. In light of Theorem~\ref{T:annealed}, it is natural to conjecture that $\rho_B(t)$ decays exactly as prescribed in Theorem~\ref{T:annealed} with $\nu=\rho_A(0)-\rho_B(0)$ and $\gamma=\infty$, whence we obtain not only the logarithmic order of decay as in \cite{BL91}, but also the constant pre-factor. However we will not address this issue here.

\subsubsection{Random Walk Among Moving Catalysts} Instead of considering $\xi$ as a field of moving traps, we may consider it as a field of moving catalysts for a system of branching random walks which we call reactants. At time $0$, a single reactant starts at the origin which undergoes branching. Independently, each reactant performs simple symmetric random walk on $\Z^d$ with jump rate $\kappa$, and undergoes binary branching with rate $|\gamma|\xi(t,x)$ when the reactant is at position $x$ at time $t$. This model was studied by Kesten and Sidoravicius in \cite{KS03}, and in the setting of the parabolic Anderson model, studied by G\"artner and den Hollander in \cite{GH06}. For the catalytic model, $\gamma$ is negative in (\ref{quenchedZ}), (\ref{annealedZ}), (\ref{FKrep}) and (\ref{annealedequiv}), and $\Z^\gamma_{t,\xi}$ and $\E^\xi[Z^\gamma_{t,\xi}]$ now represent the quenched, resp.\ annealed, expected number of reactants at time $t$. It was shown in \cite{GH06} that $\E^\xi[Z^\gamma_{t,\xi}]$ grows double exponentially fast
(i.e., $t^{-1}\log\log \E^\xi[Z^{\gamma}_{t,\xi}]$ tends to a positive limit as $t\to\infty$) for all $\gamma<0$ in dimensions $d=1$ and $2$. In $d\geq 3$, there exists
a critical $\gamma_{c,d}<0$ such that $\E^\xi[Z^\gamma_{t,\xi}]$ grows double exponentially for $\gamma<\gamma_{c,d}$, and grows exponentially (i.e., $t^{-1}\log \E^\xi[Z^{\gamma}_{t,\xi}]$ tends to a positive limit as $t\to\infty$) for all $\gamma \in (\gamma_{c,d}, 0)$. In the quenched case, however, it was shown in \cite{KS03} that $Z^{\gamma}_{t,\xi}$ only exhibits exponential growth (with $\log Z^\gamma_{t,\xi}$ shown to be of order $t$) regardless of the dimension $d\geq 1$ and the strength of interaction $\gamma<0$. Such dimension dependence bears similarities with our results for the trap model in Theorems~\ref{T:annealed} and \ref{T:quenched}.

\subsubsection{Directed Polymer in a Random Medium} We used $Z^\gamma_{t,\xi}$ to denote the survival probability, because $Z^\gamma_{t,\xi}$ and $\E^\xi[Z^\gamma_{t,\xi}]$ are in fact the quenched resp.\ annealed partition functions of a directed polymer model in a random time-dependent potential $\xi$ at inverse temperature $\gamma$. The directed polymer is modeled by $(X(s))_{0\leq s\leq t}$. In the polymer measure, a trajectory $(X(s))_{0\leq s\leq t}$ is re-weighted by the survival probability of a random walk following that trajectory in the environment $\xi$. Namely, we define a change of measure on $(X(s))_{0\leq s\leq t}$ with density $e^{-\gamma \int_0^t \xi(s,X(s))\,{\rm d}s}/Z^\gamma_{t,\xi}$
in the quenched model, and with density $\E^\xi[e^{-\gamma \int_0^t \xi(s,X(s))\,{\rm d}s}]/\E^\xi[Z^\gamma_{t,\xi}]$ in the annealed model. Qualitatively, the polymer measure favors trajectories which seek out space-time regions void of traps. However, a more quantitative geometric characterization as was carried out for the case of immobile traps (see e.g.\ \cite{S98}) is still lacking.
\vspace{1cm}

For readers interested in more background on the problem of a Brownian motion (or random walk) in time-independent potential, we refer to the book by Sznitman~\cite{S98} on Brownian motion among Poissonian obstacles,  and the survey by G\"artner and K\"onig~\cite{GK05} on the parabolic Anderson model. For readers interested in more recent studies of a random walk in time-dependent catalytic environments, we refer to the survey by G\"artner, den Hollander and Maillard~\cite{GHM09}. For readers interested in more recent studies of the trapping problem in the physics literature, we refer to the papers of Moreau, Oshanin, B\'enichou and Coppey~\cite{MOBC03, MOBC04} and the references therein.

After the completion of this paper, we learnt that the continuum analogue of our model, i.e., the study of the survival probability of a Brownian motion among a Poisson field of moving obstacles, have recently been carried out by Peres, Sinclair, Sousi, and Stauffer in~\cite{PSSS11}. See Theorems 1.1 and 3.5 therein.

\subsection{Outline}
The rest of this paper is organized as follows. Section \ref{S:annealed} is devoted to the proof of Theorem~\ref{T:annealed} on the annealed survival probability,
where the so-called Pascal principle will be introduced. In Section~\ref{S:quenched1}, we give a preliminary upper bound on the quenched survival probability in dimensions $1$ and $2$, as well as an upper bound for a semi-annealed system. Lastly, in Section~\ref{S:qexistence}, we prove the existence of the quenched Lyapunov exponent in Theorems~\ref{T:quenched} and \ref{T:quenchedPAM} via a shape theorem, and we show that the quenched Lyapunov exponent is always positive.


\section{Annealed survival probability}\label{S:annealed}
In this section, we prove Theorem~\ref{T:annealed}. We start with a proof in Section~\ref{S:lyapunov} of the existence of the annealed Lyapunov exponent $\lambda_{d,\gamma,\kappa,\rho, \nu}$. Our proof follows the same argument as for the catalytic model with $\gamma<0$ in G\"artner and den Hollander~\cite{GH06}, which is based on a special representation of $\E^\xi[Z^\gamma_{t,\xi}]$ after integrating out the Poisson random field $\xi$, which then allows us to apply the subadditivity lemma. In Section~\ref{S:annealed1}, we prove Theorem~\ref{T:annealed} for the special case $\kappa=0$, i.e., $X\equiv 0$, relying
on exact calculations. Sections~\ref{S:annealed2} and \ref{S:pascal} prove respectively the lower and upper bound on $\E^\xi[Z^\gamma_{t,\xi}]$ in Theorem~\ref{T:annealed}, for $d=1,2$ and general $\kappa>0$. The lower bound is obtained by creating a space-time box void of traps and forcing $X$ to stay inside the box, while the upper bound is based on the so-called Pascal principle, first introduced in the physics literature by Moreau et al~\cite{MOBC03, MOBC04}. In Section~\ref{S:pascal}, we will also prove the aforementioned Corollary~\ref{C:range} on the range of a symmetric random walk.

\subsection{Existence of the annealed Lyapunov exponent}\label{S:lyapunov}

In this section, we prove the existence of the annealed Lyapunov exponent
\beq\label{annlyapex}
\lambda=\lambda_{d,\gamma,\kappa,\rho,\nu}:=-\lim_{t\to\infty} \frac{1}{t}\log \E^\xi[Z^\gamma_{t,\xi}].
\eeq
\smallskip
{\bf Remark.}
 Clearly $\lambda\geq 0$, and Theorem~\ref{T:annealed} will imply that $\lambda$ always equals $0$ in dimensions $d=1,2$. For $d\geq 3$, the lower bound for
the quenched survival probability in Theorem~\ref{T:quenched} will imply that $\lambda<\gamma\nu+\kappa<\infty$, while an exact calculation of $\lambda$ for
the case $\kappa=0$ in Section~\ref{S:annealed1} and the Pascal principle in Section~\ref{S:pascal} will imply that $\lambda>0$ for all $\gamma, \nu, \rho>0$
and $\kappa\geq 0$.
\medskip

\noindent
{\bf Proof of (\ref{annlyapex}).} The proof is similar to that for the catalytic model with $\gamma<0$ in \cite{GH06}.
As in \cite{GH06}, we can integrate out the Poisson system $\xi$ to obtain
\beq \label{urep2}
\E^\xi[Z^\gamma_{t,\xi}] = \E^\xi[u(t,0)]=\E^X_0\E^\xi\!\left[\exp\!\left\{\!-\gamma\!\int_0^t \!\!\xi(t-s, X(s))\,{\rm d}s\!\right\}\! \right]\! =\E^X_0\Big[ \exp \Big\{\nu\!\! \sum_{y\in\Z^d} (v_X(t,y)-1)\Big\}\Big],
\eeq
where conditional on $X$,
\beq \label{vrep1}
v_X(t, y) = \E^Y_y\left[\exp \left\{ - \gamma\int_0^t \delta_0(Y(s)-X(t-s))\,{\rm d}s\right\}   \right]
\eeq
with $\E^Y_y[\cdot]$ denoting expectation with respect to a simple symmetric random walk $Y$ with jump rate $\rho$
and $Y(0) = y$. By the Feynman-Kac formula, $(v_X(t, y))_{t\geq 0,y\in\Z^d}$  solves the equation
\beq \label{vX}
\begin{aligned}
\frac{\partial}{\partial t} v_X(t, y) & = \rho\Delta v_X(t, y) - \gamma \delta_{X(t)}(y)\, v_X(t, y), \\
v_X(0,\cdot) & \equiv 1,
\end{aligned}
\qquad y\in\Z^d, \ t\geq 0,
\eeq
which implies that $\Sigma_X(t) := \sum_{y \in \Z^d}(v_X(t, y) -1)$ is the solution of the equation
\beq
\begin{aligned}
\frac{\rm d}{{\rm d}t}\Sigma_X(t)  & = - \gamma v_X(t, X(t)), \\
\Sigma_X(0) & = 0.
\end{aligned}
\eeq
Hence, $\Sigma_X(t) = -\gamma \int_0^t v_X(s, X(s))\, {\rm d}s$, and the representation (\ref{urep2}) becomes
\beq \label{urep3}
\E^\xi[Z^\gamma_{t,\xi}] = \E^X_0\left[ \exp \left\{-\nu\gamma \int_0^t v_X(s, X(s))\, {\rm d}s \right\}\right].
\eeq
We now observe that for $t_1, t_2>0$,
\begin{eqnarray}
\E^\xi[Z^\gamma_{t_1+t_2,\xi}] &=& \E^X_0\left[ \exp \left\{-\nu\gamma \int_0^{t_1} v_X(s, X(s))\, {\rm d}s \right\}\exp \left\{-\nu\gamma \int_{t_1}^{t_1+t_2} v_X(s, X(s))\, {\rm d}s \right\}\right] \nn\\
&\geq& \E^X_0\left[ \exp \left\{-\nu\gamma \int_0^{t_1} v_X(s, X(s))\, {\rm d}s \right\}\exp \left\{-\nu\gamma \int_{0}^{t_2} v_{\theta_{t_1}X}(s, (\theta_{t_1}X)(s))\, {\rm d}s \right\}\right] \nn\\
&=& \E^\xi[Z^\gamma_{t_1,\xi}] \E^\xi[Z^\gamma_{t_2,\xi}], \label{subadd}
\end{eqnarray}
where $\theta_{t_1}X:=((\theta_{t_1}X)(s))_{s\geq 0}=(X(t_1+s)-X(t_1))_{s\geq 0}$, we used the independence of $(X(s))_{0\leq s\leq t_1}$ and $((\theta_{t_1}X)(s))_{0\leq s\leq t_2}$, and the fact that for $s>t_1$,
\begin{eqnarray*}
v_X(s, X(s)) &=& \E^Y_{X(s)}\left[\exp \left\{ - \gamma\int_0^s \delta_0(Y(r)-X(s-r))\, {\rm d}r\right\}   \right] \\
&\leq& \E^Y_{X(s)}\left[\exp \left\{ - \gamma\int_0^{s-t_1} \delta_0(Y(r)-X(s-r))\, {\rm d}r\right\}   \right] = v_{\theta_{t_1}X}(s-t_1, (\theta_{t_1}X)(s-t_1)).
\end{eqnarray*}
From (\ref{subadd}), we deduce that $-\log\E^\xi[Z^\gamma_{t,\xi}]$ is subadditive in $t$, and hence the limit in (\ref{annlyapex}) exists and
\beq
\lambda_{d,\gamma,\kappa,\rho,\nu}=-\sup_{t>0} \frac{1}{t}\log \E^\xi[Z^\gamma_{t,\xi}].
\eeq
\qed


\subsection{Special case $\kappa = 0$}\label{S:annealed1}
In this section, we prove Theorem~\ref{T:annealed} for the case $\kappa=0$, which will be useful for lower bounding $\E^\xi[Z^\gamma_{t,\xi}]$
for general $\kappa>0$, as well as for providing an upper bound on $\E^\xi[Z^\gamma_{t,\xi}]$ by the Pascal principle.
\medskip

\noindent
{\bf Proof of Theorem~\ref{T:annealed} for $\kappa=0$.}
We first treat the case $\gamma\in (0,\infty)$. When $\kappa = 0$,
(\ref{urep3}) becomes
\beq \label{urep4}
\E^\xi[Z^\gamma_{t,\xi}] = \exp \left\{-\nu\gamma \int_0^t v_0(s, 0)\, {\rm d}s \right\},
\eeq
where $v_0$ is the solution of (\ref{vX}) with $X\equiv 0$. It then suffices to analyze the asymptotics of $v_0(t, 0)$ as $t\to \infty$.
Note that the representation (\ref{vrep1}) for $v_0(t, 0)$ becomes
\beq \label{vrep2.5}
v_0(t, 0) = \E^Y_0[e^{-\gamma \int_0^t \delta_0(Y(s))\, {\rm d}s}],
\eeq
which is the Laplace transform of the local time of $Y$ at the origin. For $d=1,2$, $v_0(t,0) \downarrow 0$ as $t\uparrow \infty$ by the
recurrence of simple random walks, while for $d\geq 3$, $v_0(t,0) \downarrow C_d$ for some $C_d>0$ by transience.

By Duhamel's principle (see e.g.~\cite[pp.~49]{E10} for a continuous-space version),
we have the following integral representation for the solution $v_X$ of (\ref{vX}),
\beq\label{vintrep}
v_X(t,y) = 1 -\gamma \int_0^t p_{\rho s}\big(y-X(t-s)\big)\, v_X\big(t-s, X(t-s)\big)\, {\rm d}s,
\eeq
where $p_s(\cdot)$ is the transition probability kernel of a rate $1$ simple symmetric random walk on $\Z^d$. When $X\equiv 0$,
we obtain
\beq \label{vrep4}
v_0(t, 0) = 1 - \gamma \int_0^t p_{\rho s}(0) v_0(t-s, 0)\,  {\rm d}s.
\eeq
Denote the Laplace transforms (in $t$) of $v_0(t, 0)$ and $p_t(0)$ by
\beq \label{vplaplace}
\hat v_0(\lambda) = \int_0^\infty e^{-\lambda t} v_0(t,0)\, {\rm d}t, \qquad \hat p(\lambda) = \int_0^\infty
e^{-\lambda t} p_t(0)\,  {\rm d}t.
\eeq
Taking Laplace transform in (\ref{vrep4}) and solving for $\hat v_0(\lambda)$ then gives
\beq
\hat v_0(\lambda) = \frac{1}{\lambda}\cdot \frac{\rho}{\rho+ \gamma\, \hat p(\lambda/\rho)}. \label{vlaplace0}
\eeq
We can apply the local central limit theorem for continuous time simple random walks in $d=1$ and 2 (i.e., $p_t(0) =
\left(\frac{d}{2\pi t}\right)^{d/2}(1+o(1))$ as $t\to\infty$) to obtain the following asymptotics for $\hat p(\lambda)$ as $\lambda \downarrow 0$,
\beq \label{plapasymp}
\hat p(\lambda)\ = \
\left\{
\begin{aligned}
\frac{1}{\sqrt{2\lambda}}(1+o(1)), & \qquad \quad d=1, \\
\frac{\ln\left(\frac{1}{\lambda}\right)}{\pi}(1+o(1)), & \qquad \quad d=2, \\
G_d(0)(1+o(1)), & \qquad \quad d\geq 3,
\end{aligned}
\right.
\eeq
with $G_d(0)=\int_0^\infty p_t(0)\, {\rm d}t$, which translates into the following asymptotics for $\hat v_0(\lambda)$ as $\lambda \downarrow 0$:
\beq \label{vlapasymp}
\hat v_0(\lambda)\ = \
\left\{
\begin{aligned}
\frac{\sqrt{2\rho}}{\gamma} \cdot \frac{1}{\sqrt \lambda}(1+o(1)),\quad & \qquad \quad d=1, \\
\frac{\pi\rho}{\gamma} \cdot \frac{1}{\lambda \ln\left(\frac{1}{\lambda}\right)}(1+o(1)), & \qquad \quad d=2, \\
\frac{\rho}{\rho+\gamma G_d(0)}\cdot \frac{1}{\lambda}(1+o(1)), & \qquad \quad d\geq 3.
\end{aligned}
\right.
\eeq
Since $v_0(t,0)$ is monotonically decreasing in $t$ by (\ref{vrep2.5}), by Karamata's Tauberian theorem (see e.g.~\cite[Chap.~XIII.5, Thm.~4]{F66}),
we have the following asymptotics for $v_0(t,0)$ as $t\to \infty$,
\beq \label{vasymp}
v_0(t, 0)\ = \
\left\{
\begin{aligned}
\frac{1}{\gamma} \sqrt{\frac{2\rho}{\pi}} \cdot \frac{1}{\sqrt t} (1+o(1)), & \qquad \quad d=1, \\
\frac{\pi\rho}{\gamma} \cdot \frac{1}{\ln t}(1+o(1)), & \qquad \quad d=2, \\
\frac{\rho}{\rho+\gamma G_d(0)}(1+o(1)), & \qquad \quad d\geq 3,
\end{aligned}
\right.
\eeq
which by (\ref{urep4}) implies Theorem~\ref{T:annealed} for $\kappa=0$ and $\gamma\in (0,\infty)$.

When $\kappa=0$ and $\gamma=\infty$, we have
$$
\E^\xi[Z^\gamma_{t,\xi}] = \P\Big(\xi(s,0)=0\ \forall\, s\in [0,t]\Big) = \exp\Big\{-\nu\sum_{y\in\Z^d}\psi(t,y)\Big\},
$$
where $\psi(t,y) = \P^Y_y(\exists\, s\in [0,t]: Y(s)=0)$ for a jump rate $\rho$ simple symmetric random walk $Y$ starting from $y$.
Note further that $\psi(t, y)$ solves the parabolic equation
\beq
\frac{\partial}{\partial t} \psi(t, y) = \rho \Delta \psi(t, y),  \qquad y\neq 0, t\geq 0,
\eeq
with boundary conditions $\psi(\cdot, 0)\equiv 1$ and $\psi(0, \cdot) \equiv 0$. Therefore $\sum_{y\in\Z^d} \psi(t, y)$ solves the equation
\beq
\frac{\rm d}{{\rm d} t} \sum_{y\in\Z^d}\psi(t, y) = -\rho\Delta \psi(t, 0)= \rho (1-\psi(t, e_1)) = \rho \phi(t, e_1),
\eeq
where $e_1=(1, 0, \cdots,0)$, $\phi(t, e_1) := 1-\psi(t, e_1)$, and we have used the fact that $\sum_{x\in\Z^d} \Delta
\psi(t, x)=0$ and the symmetry of the simple symmetric random walk. Therefore
\beq\label{gammainfty}
\E^\xi[Z^\gamma_{t,\xi}] = \exp\left\{ -\nu\rho \int_0^t \phi(s, e_1) \,{\rm d}s\right\}.
\eeq
By generating function calculations and Tauberian theorems (see e.g.~\cite[Sec.~2.4]{L96} or \cite[Sec.~32, P3]{S76}), it is known that $\phi(t, e_1)$, which is the probability that a rate 1
simple random walk starting from $e_1$ does not hit 0 before time $\rho t$, has the asymptotics $\phi(t, e_1)=\sqrt{\frac{2}{\pi\rho t}}(1+o(1))$ for
$d=1$, $\phi(t, e_1)=\frac{\pi}{\ln t}(1+o(1))$ for $d=2$, and $\phi(t,e_1)=G_d(0)^{-1}(1+o(1))$ for $d\geq 3$. Therefore as $t\to\infty$,
\beq\label{Zinftygamma}
\log \E^\xi[Z^\gamma_{t,\xi}] = \left\{
\begin{aligned}
-\nu\sqrt{\frac{8\rho t}{\pi}}(1+o(1)), \quad \quad & d=1,\\
-\nu\pi\frac{\rho t}{\ln t}(1+o(1)), \quad \quad & d=2, \\
-\nu\frac{\rho t}{G_d(0)}(1+o(1)), \quad \quad & d\geq 3,
\end{aligned}
\right.
\eeq
which proves Theorem~\ref{T:annealed} for $\kappa=0$ and $\gamma=\infty$.
\qed
\bigskip

\noindent
{\bf Remark.} When $\kappa=0$ so that $X\equiv 0$, the representation (\ref{urep4}) allows us to easily compute
the Laplace transform of $D_t:= \frac{1}{t}\int_0^t \xi(s,0)\, {\rm d}s$, since $\E^\xi[Z^\gamma_{t,\xi}]=\E^\xi[\exp\{-\gamma t D_t\}]$. By replacing $\gamma t$ with
a suitable scale $\lambda t/a_t$, where $\lambda\in\R$, $a_t=\sqrt{t}$ for $d=1$, $a_t=\log t$ for $d=2$, and $a_t=1$ for $d\geq 3$, we can identify
$$
\Psi(-\lambda) := \lim_{t\to\infty} \frac{a_t}{t} \E^\xi\big[\exp\big\{-\frac{\lambda t}{a_t}D_t\big\}\big]
$$
using the asymptotics in (\ref{vasymp}). As shown in Cox and Griffeath \cite{CG84}, applying the G\"artner-Ellis theorem then leads to a large deviation principle
for $D_t$ with scale $t/a_t$, except that in \cite{CG84}, the derivation of $\Psi(-\lambda)$ was by Taylor expansion in $\lambda$, which can be greatly simplified if we use the representation from (\ref{urep4}) instead.

\subsection{Lower bound on the annealed survival probability}\label{S:annealed2}
In this section, we prove the lower bound on $\E^\xi[Z^\gamma_{t,\xi}]$ in Theorem~\ref{T:annealed} for dimensions $d=1$ and $2$, i.e.,
\blem\label{lem:kgeq0}
For all $\gamma\in (0,\infty]$, $\kappa\geq 0$, $\rho>0$ and $\nu>0$, we have
\beq\label{lowerbd}
\begin{aligned}
\liminf_{t\to\infty} \frac{1}{\sqrt t}\log\E^\xi[Z^\gamma_{t,\xi}]\ \ &\geq\ -\nu \sqrt{\frac{8\rho}{\pi}}, &\qquad d=1, \\
\liminf_{t\to\infty} \frac{\ln t}{t}\log\E^\xi[Z^\gamma_{t,\xi}] \ &\geq\ \ -\nu\pi\rho, &\qquad d=2.
\end{aligned}
\eeq
\elem
{\bf Proof.} The basic strategy is the same as for the case of immobile traps, namely, we force the environment $\xi$ to create a ball
$B_{R_t}$ of radius $R_t$ around the origin, which remains void of traps up to time $t$, and we force the random walk $X$ to stay inside
$B_{R_t}$ up to time $t$. This leads to a lower bound on the survival probability that is independent of $\gamma\in (0,\infty]$ and
$\kappa\geq 0$. Surprisingly, in dimensions $d=1$ and $2$, this lower bound turns out to be sharp, which can be attributed to the larger
fluctuation of the random field $\xi$ in $d=1$ and $2$, which makes it easier to create space-time regions void of traps. Note that it is clearly more costly
to maintain the same space-time region void of traps than in the case when the traps are immobile.

Recall that $\xi$ is the counting field of a family of independent random walks $\{Y^y_j\}_{y\in\Z^d, 1\leq j\leq N_y}$,
where $\{N_y\}_{y\in\Z^d}$ are i.i.d.\ Poisson random variables with mean $\nu$. Let $B_r$ denote the ball of radius
$r$, i.e., $B_r=\{x\in\Z^d : \Vert x\Vert_\infty\leq r\}$. For a scale function $1<\!<R_t<\!<\sqrt{t}$ to be chosen later, let $E_t$
denote the event that $N_y=0$ for all $y\in B_{R_t}$. Let $F_t$ denote the event that $Y^y_j(s) \notin B_{R_t}$ for all
$y\notin B_{R_t},$ $1\leq j\leq N_y,$ and $s \in [0,t];$ furthermore, let $G_t$ denote the event that $X$ with $X(0)=0$ does not leave $B_{R_t}$ before time $t$.
Then by (\ref{annealedZ}),
\beq\label{422}
\E^\xi[Z^\gamma_{t,\xi}] \geq \P(E_t \cap F_t \cap G_{t}) = \P(E_t)\P(F_t)\P(G_{t}).
\eeq

Note that $\P(E_t) = e^{-\nu (2R_t+1)^d}$. To estimate $\P(G_{t})$, note that by Donsker's invariance principle, if $1<\!<R_t<\!<\sqrt{t}$ as $t\to\infty$,
then there exists $\alpha>0$ such that for all $t$ sufficiently large,
\beq
\inf_{x\in B_{\sqrt t/2}}\P\left(X(s) \in B_{\sqrt t}\ \forall\ s\in [0,t]\, ,\ X(t) \in B_{\sqrt t/2} \Big| X(0)=x\right) \geq \alpha.
\eeq
By partitioning $[0,t]$ into intervals of length $R_t^2$ and applying the Markov property, we obtain
\begin{eqnarray}
\P(G_{t}) &\geq& \P\Big(X(s)\in B_{R_t}\ \forall\ s\in [(i-1)R_t^2, iR_t^2], \text{ and } X(i R_t^2)\in B_{R_t/2},\ i=1, 2, \cdots,
\lceil t/R_t^2\rceil\Big) \nonumber\\
&\geq& \alpha^{t/R_t^2} = e^{t\ln \alpha/R_t^2}.
\end{eqnarray}

To estimate $\P(F_t)$, let $\tilde F_t$ denote the event that $Y^y_j(s) \neq 0$ for all $y\in
\Z^d$, $1\leq j\leq N_y$, and $s\in [0,t]$. Note that $\P(\tilde F_t)$ is precisely the annealed survival probability $\E^\xi[Z^\gamma_{t,\xi}]$ when $\kappa=0$ and $\gamma=\infty$, which satisfies the asymptotics in Theorem~\ref{T:annealed} by our calculations in Section \ref{S:annealed1}. We next compare $\P(F_t)$ with $\P(\tilde F_t)$.

For a jump rate $\rho$ simple random walk $Y$ starting from $y\in\Z^d$,
let $\tau_{B_{R_t}}$ denote the stopping time when $Y$ first enters $B_{R_t}$, and $\tau_0$ the stopping time when $Y$ first visits $0$.
Then standard computations yield
\beq\label{lnPFt}
\ln \P(F_t) = -\nu \sum_{y\in \Z^d\backslash B_{R_t}} \P^Y_y(\tau_{B_{R_t}}\leq t),
\eeq
and a similar identity holds for $\ln\P(\tilde F_t)$ with $B_{R_t}$ replaced by $B_0$. Note that
\beq\nonumber
\sum_{y\in \Z^d\backslash B_{R_t}} \P^Y_y(\tau_{B_{R_t}}\leq t) \geq \sum_{y\in \Z^d\backslash B_{R_t}}
\P^Y_y(\tau_0\leq t) = \sum_{y\in\Z^d} \P^Y_y(\tau_0\leq t) - \sum_{y\in B_{R_t}}\P^Y_y(\tau_0\leq t).
\eeq
Hence
\beq\label{431.1}
\ln \P(F_t) \leq \ln \P(\tilde F_t) + \nu \sum_{y\in B_{R_t}} \P^Y_y(\tau_0\leq t)
\leq \ln \P(\tilde F_t) + \nu (2R_t+1)^d.
\eeq
On the other hand, for $\eps>0$, we have
$$
\sum_{y\in\Z^d} \P^Y_y(\tau_0\leq t+\eps t) \geq \!\!\!\! \sum_{y\in\Z^d\backslash B_{R_t}}
\P^Y_y\big(\tau_{B_{R_t}}\leq t, \tau_0 \leq t+\eps t\big)
\geq \inf_{z\in \partial B_{R_t}}\P^Y_z(\tau_0\leq \eps t)
\!\!\!\!\sum_{y\in\Z^d\backslash B_{R_t}} \P^Y_y\big(\tau_{B_{R_t}}\leq t),
$$
where we used the strong Markov property. Therefore
$$
\sum_{y\in \Z^d\backslash B_{R_t}} \P^Y_y(\tau_{B_{R_t}}\leq t) \leq \frac{\sum_{y\in\Z^d}\P^Y_y(\tau_0\leq t+\eps
  t)} {\inf_{z\in \partial B_{R_t}}\P^Y_z(\tau_0\leq \eps t)},
$$
and hence by (\ref{lnPFt}),
\beq\label{431.2}
\ln \P(F_t) \geq \frac{\ln \P(\tilde F_{t+\eps t})}{\inf_{z\in \partial B_{R_t}} \P^Y_z(\tau_0\leq \eps t)}.
\eeq

We now choose $R_t$ for $d=1$ and $2$. For $d=1$, let $R_t= \sqrt{t/\ln t}$, which is by no means the
unique scale appropriate. Clearly $\inf_{z\in \partial B_{\sqrt{t/\ln t}}}\P^Y_z(\tau_0\leq \eps t) \to 1$ as
$t\to\infty$. By (\ref{431.1})--(\ref{431.2}), the fact that $\P(\tilde F_t)$ satisfies the asymptotics in Theorem~\ref{T:annealed}
for $\kappa=0$ and $\gamma=\infty$, and that $\eps>0$ can be made arbitrarily small, we obtain
$$
\ln \P(F_t) = -\nu \sqrt{\frac{8\rho t}{\pi}} (1+o(1)) = \ln \P(\tilde F_t).
$$
Furthermore, for $R_t = \sqrt{t/\ln t}$ we have
$$
\ln \P(E_t) =-\nu (2\sqrt{t/\ln t}+1) \qquad \mbox{and} \qquad \ln \P(G_{t}) \geq \ln \alpha \ln t,
$$
whence substituting these asymptotics into (\ref{422}) gives (\ref{lowerbd}) for $d=1$.

For $d=2$, let $R_t = \ln t$. Then we have $\inf_{z\in \partial B_{\ln t}}\P^Y_z(\tau_0\leq \eps t) \to 1$ as
$t\to\infty$, which is an easy consequence of~\cite[Exercise 1.6.8]{L96}. By the same
argument as for $d=1$, we have
$$
\ln \P(F_t) =-\nu\pi\rho \frac{t}{\ln t} (1+o(1))= \ln \P(\tilde F_t).
$$
Together with the asymptotics
$$
\ln \P(E_t) = -\nu (2\ln t+1)^2 \qquad \mbox{and} \qquad \ln \P(G_{t}) \geq \frac{t \ln\alpha}{\ln^2t},
$$
we deduce from (\ref{422}) the desired bound in (\ref{lowerbd}) for $d=2$.
\qed
\bigskip

\subsection{Upper bound on the annealed surivival probability: the Pascal principle}\label{S:pascal}
In this section, we present an upper bound on the annealed survival probability, called the Pascal principle.

\bprop{\bf [Pascal principle]}\label{P:pascal} Let $\xi$ be the random field generated by a collection of irreducible symmetric random walks
$\{Y^y_j\}_{y\in\Z^d, 1\leq j\leq N_y}$ on $\Z^d$ with jump rate $\rho>0$. Then for all piecewise constant $X:[0,t]\to \Z^d$ with a finite number of discontinuities,
we have
\beq\label{Pascal0}
\E^\xi\left[\exp\left\{-\gamma\int_0^t \xi(s,X(s))\,{\rm d}s\right\} \right]
\leq \E^\xi\left[\exp\left\{-\gamma\int_0^t \xi(s,0)\,{\rm d}s\right\}\right].
\eeq
\eprop
In words, conditional on the random walk $X$, the annealed survival probability is maximized when $X\equiv 0$. The discrete time version of this result was first proved by Moreau et al in \cite{MOBC03, MOBC04}, where they named it the {\em Pascal principle}, because Pascal once asserted that {\it all misfortune of men comes from the fact that he does not stay peacefully in his room}. The Pascal principle together with the
proof of Theorem~\ref{T:annealed} for $\kappa=0$ in Section~\ref{S:annealed1} imply the desired upper bound on the annealed survival probability in Theorem~\ref{T:annealed} for dimensions $d=1,2$, and it also shows that for $d\geq 3$, the annealed Lyapunov exponent $\lambda_{d,\gamma,\kappa,\rho,\nu}$ is always bounded from below by $\lambda_{d,\gamma,0,\rho,\nu}=\nu\gamma/(1+\frac{\gamma G_d(0)}{\rho})$.

We present below the proof of the discrete time version of the Pascal principle from \cite{MOBC04}, which being written as a physics paper, can be hard for the reader to separate the rigorous arguments from the non-rigorous ones. We then deduce the continuous time
version, Proposition~\ref{P:pascal}, by discrete approximation. As a byproduct, we will show in Corollary \ref{C:range} that the expected cardinality of the range of a continuous time symmetric random walk
increases under perturbation by a deterministic path.

Moreau et al considered in \cite{MOBC04} a discrete time random walk among a Poisson field of moving traps, defined as follows. Let $\bar X$ be a discrete time mean zero random walk on $\Z^d$ with $\bar X_0=0$. Let $\{N_y\}_{y\in\Z^d}$ be i.i.d.\ Poisson random variables with mean $\nu$, and let $\{\bar Y^y_j\}_{y\in\Z^d, 1\leq j\leq N_y}$  be a family of independent {\it symmetric} random walks on $\Z^d$ where $\bar Y^y_j$ denotes the $j$-th random walk starting from $y$ at time $0$. Let
\beq
\bar\xi(n,x) := \sum_{y\in \Z^d, 1\leq j\leq N_y} \delta_x(\bar Y^y_j(n)).
\eeq
Fix $0\leq q\leq 1$, which will be the trapping probability. The dynamics of $\bar X$ is such that $\bar X$
moves independently of the traps $\{\bar Y^y_j\}_{y\in\Z^d, 1\leq j\leq N_y}$, and at each time $n\geq 0$, $\bar X$ is killed
with probability $1-(1-q)^{\bar\xi(n, \bar X(n))}$. Namely, each trap at the time-space lattice site $(n,\bar X(n))$
tries independently to capture $\bar X$ with probability $q$. Given a realization of $\bar X$, let
$\bar\sigma^{\bar X}(n)$ denote the probability that $\bar X$ has survived till time $n$. Then analogous to (\ref{urep2}),
we have
\beq
\bar\sigma^{\bar X}(n) = \E^{\bar\xi}\left[(1-q)^{\sum_{i=0}^n \bar\xi(i, \bar X(i))}\right]
= \exp\Big\{-\nu \sum_{y\in\Z^d} \bar w^{q,\bar X}(n,y)\Big\},
\eeq
where if we let $\bar Y$ denote a random walk with the same jump kernel as $\bar Y^y_j$, then
\beq\label{wbar}
\bar w^{q, \bar X}(n,y) := 1- \E^{\bar Y}_y\Big[(1-q)^{\sum_{i=0}^n 1_{\{\bar Y(i)=\bar X(i)\}}}\Big].
\eeq
The main result we need from Moreau et al \cite{MOBC04} is the following discrete time Pascal principle.
\blem{\bf [Pascal principle in discrete time \cite{MOBC04}]}\label{L:pascal1}\\
Let $\bar Y$ be an irreducible symmetric random walk on $\Z^d$ with $\P^{\bar Y}_0(\bar Y(1)=0) \geq 1/2$. Then
for all $q\in [0,1]$, $n\in \N_0$ and $\bar X : \N_0\to \Z^d$, we have
\beq\label{pascal1}
\sum_{y\in\Z^d} \bar w^{q,\bar X}(n,y) \geq \sum_{y\in\Z^d} \bar w^{q,0}(n,y),
\eeq
and hence $\bar\sigma^{\bar X}(n) \leq \bar\sigma^0(n)$, where $\bar w^{q,0}$ and $\bar\sigma^0$ denote $\bar w^{q,\bar X}$ and $\bar\sigma^{\bar X}$ with $\bar X\equiv 0$.
\elem
{\bf Proof.} The argument we present here is extracted from \cite{MOBC04}. First note that the assumption
$\bar Y$ is symmetric implies that the Fourier transform $f(k):=\E^{\bar Y}_0[e^{i\langle k, \bar Y(1)\rangle}]$ is real for all $k\in [-\pi,\pi]^d$.
The assumption $\P^{\bar Y}_0(\bar Y(1)=0)\geq 1/2$ guarantees that $f(k)\in [0,1]$. If we let $p^{\bar Y}_n(y)$ denote the $n$-step transition probability kernel of
$\bar Y$, then by Fourier inversion, we have
\beq\label{monotone}
\begin{aligned}
p^{\bar Y}_n(0) \   &\geq\ p^{\bar Y}_n(y), \\
p^{\bar Y}_n(0) \ \ &\geq\ p^{\bar Y}_{n+1}(0)
\end{aligned}
\qquad \mbox{for all } n\geq 0,\ y\in\Z^d.
\eeq

If we now regard $\bar X$ as a trap, then $\bar w^{q,\bar X}(n,y)$ can be interpreted as the
probability that a random walk $\bar Y$ starting from $y$ gets trapped by $\bar X$ by time $n$, where
each time $\bar Y$ and $\bar X$ coincide, $\bar Y$ is trapped by $\bar X$ with probability $q$. More precisely,
let $Z_i$, $i\in\N_0$, be i.i.d.\ Bernoulli random variables with mean $q$, where $Z_i=1$ means that the trap at $(i, \bar X(i))$ is open. Then $\bar X$ is killed at the stopping time
\beq
\tau_{\bar X}(\bar Y) := \min\{i\geq 0: \bar Y(i)=\bar X(i), Z_i=1\},
\eeq
and $\bar w^{q,\bar X}(n,y) = \P^{\bar Y}_y(\tau_{\bar X} \leq n)$.

We examine the following auxiliary quantity, where by decomposition with respect to $\tau_{\bar X}$, we have
\beq\label{decx}
\begin{aligned}
q= \sum_{y\in\Z^d} \P^{\bar Y}_y\Big(\bar Y(n)=\bar X(n), Z_n=1\Big)
&=\ \sum_{k=0}^{n-1} \sum_{y\in\Z^d} \P^{\bar Y}_y(\tau_{\bar X}=k) p^{\bar Y}_{n-k}(\bar X(n)-\bar X(k))\, q
+ \sum_{y\in\Z^d} \P^{\bar Y}_y(\tau_{\bar X}=n) \\
&\leq\ q \sum_{k=0}^{n-1} \sum_{y\in\Z^d} \P^{\bar Y}_y(\tau_{\bar X}=k) p^{\bar Y}_{n-k}(0)
+ \sum_{y\in\Z^d} \P^{\bar Y}_y(\tau_{\bar X}=n),
\end{aligned}
\eeq
where in the inequality we used (\ref{monotone}). Similarly, when $\bar X$ is replaced by $\bar X\equiv 0$,
we have
\beq\label{dec0}
q=\sum_{y\in\Z^d} \P^{\bar Y}_y\Big(\bar Y(n)=0, Z_n=1\Big)
= q \sum_{k=0}^{n-1} \sum_{y\in\Z^d} \P^{\bar Y}_y(\tau_{0}=k) p^{\bar Y}_{n-k}(0)
+ \sum_{y\in\Z^d} \P^{\bar Y}_y(\tau_{0}=n).
\eeq
Denote
\beq
\begin{aligned}
S^{\bar X}_n\ &:=\ \sum_{y\in\Z^d} \bar w^{q, \bar X}(n,y) = \sum_{y\in\Z^d} \P^{\bar Y}_y(\tau_{\bar X}\leq n), \\
S^{0}_n     \ &:=\ \sum_{y\in\Z^d} \bar w^{q, 0}(n,y) = \sum_{y\in\Z^d} \P^{\bar Y}_y(\tau_{0}\leq n).
\end{aligned}
\eeq
Note that $S^{\bar X}_0=S^{0}_0=q$, and $\sum_{y\in\Z^d}\P^{\bar Y}_y(\tau_{\bar X}=k)= S^{\bar X}_k-S^{\bar X}_{k-1}$,
$\sum_{y\in\Z^d}\P^{\bar Y}_y(\tau_0=k)= S^0_k-S^0_{k-1}$, where we set $S^{\bar X}_{-1}=S^0_{-1}=0$. Together with
(\ref{decx}) and (\ref{dec0}), this gives
$$
q \sum_{k=0}^{n-1} p^{\bar Y}_{n-k}(0) (S^0_k-S^0_{k-1}) + S^0_{n}-S^0_{n-1}
\leq
q \sum_{k=0}^{n-1} p^{\bar Y}_{n-k}(0) (S^{\bar X}_k-S^{\bar X}_{k-1}) + S^{\bar X}_n-S^{\bar X}_{n-1}.
$$
Rearranging terms, we obtain
\beq\label{induction}
S^{\bar X}_n-S^0_n \geq \big(1- qp^{\bar Y}_1(0)\big)(S^{\bar X}_{n-1}-S^0_{n-1})
+ q\sum_{k=0}^{n-2} \Big(p^{\bar Y}_{n-k-1}(0)-p^{\bar Y}_{n-k}(0)\Big) (S^{\bar X}_k -S^0_{k}).
\eeq
This sets up an induction bound for $S^{\bar X}_n-S^0_n$. Since $S^{\bar X}_0-S^0_0=0$, $1-q p^{\bar Y}_1(0)\geq 0$,
and $p^{\bar Y}_k(0)$ is decreasing in $k$ by (\ref{monotone}), it follows that $S^{\bar X}_n\geq S^0_n$ for
all $n\in\N_0$, which is precisely (\ref{pascal1}).
\qed
\bigskip

\noindent
{\bf Proof of Proposition~\ref{P:pascal}.} Integrating out $\xi$ on both sides of (\ref{Pascal0}) as in (\ref{urep2}) shows that
(\ref{Pascal0}) is equivalent to
\beq\label{pascal2}
\sum_{y\in\Z^d} w^{\gamma,X}(t,y) \geq \sum_{y\in\Z^d} w^{\gamma,0}(t,y),
\eeq
where
\beq
w^{\gamma,X}(t,y) := 1 - \E^Y_y\left[\exp\Big\{-\gamma \int_0^t \delta_0(Y(s)-X(s))\, {\rm d}s\Big\}\right].
\eeq
For $n\in\N$, let
$Y^{(n)}(k) = Y(\frac{k t}{n})$ and $X^{(n)}(k)= X(\frac{ t}{n})$ for $k\in\N_0$. Clearly $Y^{(n)}$ is symmetric, and for $n$ sufficiently
large, $\P^{Y^{(n)}}_0(Y^{(n)}(1)=0)\geq 1/2$. Therefore we can apply Lemma \ref{L:pascal1} with $\bar Y=Y^{(n)}$,
$\bar X=X^{(n)}$ and $q=q^{(n)} = \gamma t/n$ to obtain
\beq\label{approxineq}
\sum_{y\in\Z^d} \bar w^{\gamma t/n, X^{(n)}}(n, y) \geq \sum_{y\in\Z^d} \bar w^{\gamma t/n, 0}(n, y).
\eeq
By (\ref{wbar}) and the definition of $Y^{(n)}$ and $X^{(n)}$, we have
$$
\bar w^{\gamma t/n, X^{(n)}}(n,y) = 1-\E^{Y^{(n)}}_y\left[\Big(1-\frac{\gamma t}{n}\Big)^{\sum_{k=0}^n
    1_{\{Y^{(n)}(k)=X^{(n)}(k)\}}}\right]=1- \E^Y_y\left[\Big(1-\frac{\gamma t}{n}\Big)^{\sum_{k=0}^n 1_{\{Y(kt/n)=X(kt/n)\}}}\right].
$$
By the assumption that $X$ is a random walk path which is necessarily piecewise constant with a finite number of
discontinuities, for a.s.\ all realization of $Y$, we have
$$
\lim_{n\to\infty} \Big(1-\frac{\gamma t}{n}\Big)^{\sum_{k=0}^n 1_{\{Y(kt/n)=X(kt/n)\}}}
= \exp\left\{-\gamma \int_0^t \delta_0(Y(s)-X(s))\, {\rm d}s\right\}.
$$
Therefore by the bounded convergence theorem, $\lim_{n\to\infty} \bar w^{\gamma t/n, X^{(n)}}(n,y) = w^{\gamma,X}(t,y)$.
By the same argument, $\lim_{n\to\infty} \bar w^{\gamma t/n, 0}(n,y) = w^{\gamma,0}(t,y)$. Next we note
that $w^{\gamma t/n, X^{(n)}}(n,y)$ is the probability that $Y^{(n)}$ is trapped by $X^{(n)}$ before time $n$.
Since $Y^{(n)}$ and $X^{(n)}$ are embedded in $Y$ and $X$, we have
$w^{\gamma t/n, X^{(n)}}(n,y) \leq \P^Y_y(\tau_X\leq t)$ uniformly in $n$, where $\tau_X=\inf\{s\geq 0: Y(s)=X(s)\}$.
Clearly $\sum_{y\in\Z^d} \P^Y_y(\tau_X\leq t)<\infty$. Similarly $w^{\gamma t/n, 0}(n,y) \leq \P^Y_y(\tau_0\leq t)$ uniformly
in $n$ and $\sum_{y\in\Z^d}\P^Y_y(\tau_0\leq t)<\infty$. Therefore we can send $n\to\infty$ and apply the dominated convergence
theorem in (\ref{approxineq}), from which (\ref{pascal2}) then follows.
\qed
\bigskip

The Pascal principle in Lemma \ref{L:pascal1} and Proposition \ref{P:pascal} have the following interesting consequence for the range of
a symmetric random walk, which we denote by $R_t(X) = \{y\in \Z^d: X(s)=y \mbox{ for some } 0\leq s\leq t\}$.
\bcor{\bf [Increase of expected cardinality of range under perturbation]}\label{C:range}\\
Let $\bar Y$ and $\bar X$ be discrete time random walks as in Lemma \ref{L:pascal1}. Let $Y$ be a continuous time irreducible symmetric
random walk on $\Z^d$ with jump rate $\rho>0$, and let $X: [0,t]\to \Z^d$ be piecewise constant with a finite number of discontinuities. Then for
all $n\in \N_0$, respectively\ $t\geq 0$, we have
\beq\label{rangeineq}
\begin{aligned}
\E^{\bar Y}_0\big[|R_n(\bar Y-\bar X)|\big]\ &\geq\ \E^{\bar Y}_0\big[ |R_n(Y)|\big], \\
\E^{Y}_0\big[|R_t(Y-X)|\big]\ &\geq\ \E^{Y}_0\big[ |R_t(Y)|\big],
\end{aligned}
\eeq
where $|\cdot|$ denotes the cardinality of the set.
\ecor
{\bf Proof.} The first inequality in (\ref{rangeineq}) for discrete time random walks follows from the observation that
\beq
\begin{aligned}
\sum_{y\in \Z^d} \P^{\bar Y}_y(\tau_{\bar X}\leq n) \ &=\
\sum_{y\in\Z^d} \P^{\bar Y}_0\big(\bar Y(i)-\bar X(i)=y \mbox{\ for some\ } 0\leq i\leq n\big) = \E^{\bar Y}_0\big[|R_n(\bar
Y-\bar X)|\big], \\
\sum_{y\in \Z^d} \P^{\bar Y}_y(\tau_{0}\leq n) \ &=\
\sum_{y\in\Z^d} \P^{\bar Y}_0\big(\bar Y(i) =y \mbox{\ for some\ } 0\leq i\leq n\big) = \E^{\bar Y}_0\big[|R_n(\bar
Y)|\big],
\end{aligned}
\eeq
where $\tau_{\bar X} = \min\{i\geq 0: \bar Y_i=\bar X_i\}$ and $\tau_0=\min\{i\geq 0:\bar Y_i=0\}$, which combined with
Lemma \ref{L:pascal1} for $q=1$ gives precisely
\beq
\sum_{y\in\Z^d} \P^{\bar Y}_y(\tau_{\bar X} \leq n) \geq \sum_{y\in\Z^d} \P^{\bar Y}_y(\tau_0 \leq n).
\eeq
The continuous time case follows by similar considerations, where we apply Proposition~\ref{P:pascal} with
$\gamma =\infty$, or rather $\gamma>0$ with $\gamma\uparrow \infty$.
\qed

\section{Quenched and semi-annealed upper bounds}\label{S:quenched1}
In this section, we prove sub-exponential upper bounds on the quenched survival probability in dimensions $1$ and $2$ (the exponential upper bound in dimensions $3$ and higher follows trivially from the annealed upper bound by Jensen's inequality and Borel-Cantelli). Although they will be superseded later by a proof of exponential decay using sophisticated results of Kesten and Sidoravicius~\cite{KS05}, the proof we present here is relatively simple and self-contained. Along the way, we will also prove an upper bound (Proposition~\ref{P:panneal}) on the annealed survival probability of a random walk in a random field of traps $\xi$ with deterministic initial condition, which we call a semi-annealed bound.

\bprop\label{P:subexp}{\bf [Sub-exponential upper bound on $Z^\gamma_{t,\xi}$]} There exist constants $C_1, C_2>0$ depending on $\gamma, \kappa, \rho, \nu>0$ such that
a.s.\ with respect to $\xi$, we have
\beq\label{subexp}
\begin{aligned}
\limsup_{t\to\infty} \frac{\log t}{t} \log Z^\gamma_{t,\xi} \ \leq\ -C_1, \qquad & \qquad d=1, \\
\limsup_{t\to\infty} \frac{\log\log t}{t} \log Z^\gamma_{t,\xi}\ \leq\ -C_2, \qquad & \qquad d=2.
\end{aligned}
\eeq
The same bounds hold if we replace $Z^\gamma_{t,\xi}$ by $u(t,0)$ as in Theorem~\ref{T:quenchedPAM}.
\eprop
{\bf Proof.} The proof is based on coarse graining combined with the annealed bound in Theorem~\ref{T:annealed}. Let us focus on dimension $d=1$ first. Let $X$ be a random walk as in (\ref{quenchedZ}),
and let $M(t):=\sup_{0\leq s\leq t} |X(s)|_\infty$. The first step is to note that by basic large deviation estimates for $X$,
$$
\E^X_0\Big[\exp\Big\{-\gamma\int_0^t\xi(s,X(s))\, {\rm d}s\Big\} 1_{\{M_t\geq t\}}\Big] \leq \P^X_0(M_t\geq t) \leq e^{-Ct}
$$
for some $C>0$ depending only on $\kappa$. Therefore to show (\ref{subexp}), it suffices to prove that
\beq\label{coarse1}
\E^X_0\Big[\exp\Big\{-\gamma\int_0^t\xi(s,X(s))\, {\rm d}s\Big\} 1_{\{M_t < t\}}\Big] \leq e^{-\frac{C t}{\log t}}
\eeq
for some $C>0$ for all $t$ sufficiently large. Since the
integrand in the definition of $Z_{t,\xi}^\gamma$ is monotone in $t,$
we may even restrict our attention to $t\in\N$.

The second step is to introduce a coarse graining scale $L_t:=A\log t$ for some $A>0$, and partition the space-time region $[-2t, 2t]\times [0,t]$ into blocks of the form $\Lambda_{i,k}:=[(i-1)L_t, iL_t)\times [(k-1)L_t^2, kL_t^2)$ for $i,k\in\Z$ with $-\frac{2t}{L_t}+1\leq i\leq \frac{2t}{L_t}$ and $1\leq k\leq \frac{t}{L_t^2}$. We say a block
$$
\Lambda_{i,k} \ \mbox{is good if } \ \sum_{(i-1)L_t\leq x< iL_t} \xi((k-1)L_t^2, x)\geq \frac{\nu L_t}{2}.
$$
Since for each $s\geq 0$, $(\xi(s,x))_{x\in\Z}$ are i.i.d.\ Poisson distributed with mean $\nu$, by basic large deviation estimates for Poisson random variables,
there exists $C>0$ such that for all $t>1$,
$$
\P(\Lambda_{i,k} \mbox{ is bad}) \leq e^{-C\nu L_t}.
$$
Let $G_t(\xi)$ be the event that all the blocks $\Lambda_{i,k}$ in $[-2t, 2t]\times [0,t]$ are good. Then
$$
\P(G^c_t(\xi))\leq \frac{4t^2}{L_t^3} e^{-C\nu L_t}    =\frac{4}{A^3(\log t)^3 t^{C\nu A-2}},
$$
which is summable in $t\geq 2$, $t\in\N$, if $A$ is chosen sufficiently large. Therefore by Borel-Cantelli, a.s.\ with respect to $\xi$, for all $t\in\N$
sufficiently large, the event $G_t(\xi)$ occurs. To prove (\ref{subexp}), it then suffices to prove
\beq\label{coarse2}
1_{G_t(\xi)}\E^X_0\Big[\exp\Big\{-\gamma\int_0^t\xi(s,X(s))\, {\rm d}s\Big\} 1_{\{M_t < t\}}\Big] \leq e^{-\frac{C t}{\log t}}
\eeq
almost surely for all $t\in\N$ sufficiently large.

The third step is applying an annealing bound. More precisely, to show (\ref{coarse2}), it suffices to average over $\xi$ and show that
\beq\label{coarse3}
\E^\xi\E^X_0\Big[\exp\Big\{-\gamma\int_0^t\xi(s,X(s))\, {\rm d}s\Big\} 1_{\{M_t < t\}}1_{G_t(\xi)}\Big] \leq e^{-\frac{2C t}{\log t}}
\eeq
for some $C>0$ for all $t\in\N$ sufficiently large. Indeed, (\ref{coarse3}) implies that
$$
\P^\xi\Big( 1_{G_t(\xi)}\E^X_0\Big[\exp\Big\{-\gamma\int_0^t\xi(s,X(s))\, {\rm d}s\Big\} 1_{\{M_t < t\}}\Big] > e^{-\frac{C t}{\log t}} \Big)\leq e^{-\frac{C t}{\log t}},
$$
from which (\ref{coarse2}) then follows by Borel-Cantelli.

To prove (\ref{coarse3}), let us denote $Z_k:=\exp\big\{-\gamma\int_{(k-1)L_t^2}^{kL_t^2}\xi(s,X(s))\, {\rm d}s \big\}$, and let $\Fi_k$ be the $\sigma$-field generated by
$(X_s, \xi(s,\cdot))_{0\leq s\leq kL_t^2}$. Replacing $L_t^2$ by $t/\lfloor t/L_t^2\rfloor$ if necessary, we may assume without loss of generality that
$t/L_t^2=t/(A\log t)^2\in\N$. Then
\begin{eqnarray}
&& \E^\xi\E^X_0\Big[\exp\Big\{-\gamma\int_0^t\xi(s,X(s))\, {\rm d}s\Big\} 1_{\{M_t < t\}}1_{G_t(\xi)}\Big] = \E^\xi\E^X_0\Big[1_{\{M_t < t\}}1_{G_t(\xi)}\prod_{k=1}^{t/L_t^2} Z_k  \Big] \nn \\
&=& \E^\xi\E^X_0\Big[1_{\{M_t < t\}}1_{G_t(\xi)}\prod_{k=1}^{t/L_t^2} \frac{Z_k}{\E^\xi\E^X_0[Z_k|\Fi_{k-1}]} \prod_{k=1}^{t/L_t^2}\E^\xi\E^X_0[Z_k|\Fi_{k-1}]\Big]. \label{coarse4}
\end{eqnarray}
By Proposition \ref{P:panneal} below, on the event $|X((k-1)L_t^2)|_\infty <t$
and $\Lambda_{i,k}$ is good for all $-\frac{2t}{L_t}+1\leq i\leq \frac{2t}{L_t}$, which is an event in $\Fi_{k-1}$, we have
\beq\label{coarse5}
\E^\xi\E^X_0[Z_k|\Fi_{k-1}]= \E^\xi\E^X_0\Big[ \exp\Big\{-\gamma\int_{(k-1)L_t^2}^{kL_t^2}\xi(s,X(s))\,{\rm d}s\Big\} \Big |\Fi_{k-1}\Big] \leq e^{-C L_t}
\eeq
for some $C>0$ depending on $\gamma, \kappa, \rho, \nu$. Substituting this bound into (\ref{coarse4}) for $1\leq k\leq t/L_t^2$ and using the fact
that $\prod_{k=1}^{t/L_t^2} \frac{Z_k}{\E^\xi\E^X_0[Z_k|\Fi_{k-1}]}$ is a martingale then gives the desired bound $e^{-Ct/L_t}=e^{-Ct/(A\log t)}$
for (\ref{coarse3}).

For dimension $d=2$, the proof is similar. We choose $L_t=A\log t$ with $A$ sufficiently large. We partition the space-time region $[-2t, 2t]^2\times [0,t]$
into blocks of the form $\Lambda_{i,j,k}:=[(i-1)L_t, iL_t)\times[(j-1)L_t, jL_t)\times [(k-1)L_t^2, kL_t^2)$, and we define good blocks and bad blocks as before.
Applying Proposition \ref{P:panneal} below then gives an upper bound of $\exp\big\{-C \frac{t}{L_t^2}\frac{L^2_t}{\log L_t}\big\}= \exp\{-\frac{Ct}{\log A+\log\log t}\}$, analogous to (\ref{coarse3}).

Lastly we note that the arguments also apply to the solution of the parabolic Anderson model
$$
u(t,0)=\E^X_0\Big[\exp\Big\{-\gamma\int_0^t\xi(t-s,X(s))\, {\rm d}s\Big\}\Big].
$$
The only difference lies in passing the result (\ref{subexp}) from $t\in\N$ to $t\in\R$, due to the lack of monotonicity of $u(t,0)$ in $t$. This can be easily overcome
by the observation that for $n-1<t<n$ with $n\in\N$,
$$
u(n,0)\geq e^{-\kappa (n-t)} e^{-\gamma \int_{t}^n \xi(r,0)\, {\rm d}r} u(t,0),
$$
and the fact that almost surely $\int_{i}^{i+1} \xi(r,0)\, {\rm d}r \leq \sqrt{i}$ for all $i$ large by Borel-Cantelli because
$\int_0^1\xi(r,0)\,{\rm d}r$ has finite exponential moments.
\qed
\bigskip

The following is a partial analogue of Theorem~\ref{T:annealed} for $\xi$ with deterministic initial conditions.
\bprop{\bf [Semi-annealed upper bound]}\label{P:panneal}
Let $\xi$ be defined as in (\ref{xi}) with deterministic initial condition $(\xi(0,x))_{x\in\Z^d}$. For $L>0$ and $\vec i=(i_1,\cdots, i_d)\in\Z^d$, let $B_{L,\vec i}:=[(i_1-1)L,i_1L)\times\cdots\times [(i_d-1)L,i_dL) $. Assume that there exist $a>2$ and $\nu>0$ such that for all $\vec i\in [-3L^a, 3L^a]^d$,
$\sum_{x\in B_{L,\vec i}} \xi(0,x)\geq \nu L^d$. Then there exist constants $C_d>0$, $d\geq 1$, such that for all $L$ sufficiently large and for all $x\in\Z^d$ with $|x|_\infty\leq L$, we have
\beq\label{pann}
\E^\xi\E^X_x\Big[\exp\Big\{-\gamma\int_0^{L^2} \xi(s, X(s))\, {\rm d}s\Big\}\Big]\ \leq
\left\{
\begin{aligned}
e^{-C_1L}, \qquad &\qquad d=1, \\
e^{-C_2\frac{L^2}{\log L}}, \qquad &\qquad d=2, \\
e^{-C_d L^2}, \qquad &\qquad d\geq 3.
\end{aligned}
\right.
\eeq
The same is true if we replace $\int_0^{L^2} \xi(s,X(s))\, {\rm d}s$ by $\int_0^{L^2} \xi(s,X(L^2-s))\, {\rm d}s$.
\eprop
{\bf Proof.} The basic strategy is to dominate $(\xi(L^2/2,x))_{|x|_\infty<2L^a}$ from below by i.i.d.\ Poisson random variables, which then allows us to apply Theorem~\ref{T:annealed}. We proceed as follows.

Let $\xi$ be generated by independent random walks $(Y^y_j)_{y\in\Z^d, 1\leq j\leq \xi(0,y)}$ as in (\ref{xi}),
and let $\bar\xi$ be generated by a separate system of independent random walks $(\bar Y^y_j)_{y\in\Z^d, 1\leq j\leq \bar\xi(0,y)}$, where $(\bar\xi(0,y))_{y\in\Z^d}$
are i.i.d.\ Poisson distributed with mean $\bar\nu$. Choose any $\bar\nu\in (0,\nu)$. Then by large deviation estimates for Poisson random variables,
\begin{eqnarray}
\P^{\bar \xi}(G_L^c)&:=& \P^{\bar\xi}\Big(\sum_{x\in B_{L,\vec i}} \bar\xi(0,x) \geq \sum_{x\in B_{L,\vec i}} \xi(0,x) \mbox{ for some } \vec i\in [-3L^a, 3L^a]^d   \Big) \label{poisdom}\\
&\leq& \sum_{\vec i\in [-3L^a, 3L^a]^d}\P^{\bar \xi}\Big(\sum_{x\in B_{L,\vec i}} \bar\xi(0,x) \geq \nu L^d \Big) \leq 6^d L^{ad} e^{-C_{\nu,\bar\nu} L^d}. \nn
\end{eqnarray}
On the event $G_L$, we will construct a coupling between $(Y^y_j)_{y\in\Z^d, 1\leq j\leq \xi(0,y)}$ and $(\bar Y^y_j)_{y\in\Z^d, 1\leq j\leq \bar\xi(0,y)}$ as follows. For each walk $\bar Y^y_j$ with $1\leq j\leq \bar\xi(0,y)$ and $y\in B_{L,\vec i}$ for some $\vec i\in [-3L^a, 3L^a]^d$, we can match $\bar Y^y_j$ with a distinct walk $Y^z_k$ for some $z\in B_{L,\vec i}$ and $1\leq k\leq \xi(0,z)$, which is possible on the event $G_L$.

Independently for each pair of walks $(\bar Y^y_j, Y^z_k)$, we will couple their coordinates as follows:
For $1 \leq i \leq d,$ the $i$-th coordinates of the two walks evolve independently until the first time that their
difference is of even parity. Note that this is the case either at time $0$ already or at the first time when one of the coordinates
changes. From then on the $i$-th coordinates are coupled in such a way that they always jump at the same time
and their jumps are always opposites until the first time when the two
coordinates coincide.
From that time onward the two coordinates always perform the same jumps at the same time.
For walks in the $\xi$ and $\bar\xi$ system which have not been paired up, we let them evolve independently. Note that such a coupling preserves the law of $\xi$ (resp.\ $\bar\xi$), and each coupled pair $(\bar Y^y_j, Y^z_k)$ is successfully coupled in the sense that $\bar Y^y_j(L^2/2)=Y^z_k(L^2/2)$ if the trajectory of $\bar Y^y_j$ is in the event
$$
E^y_j :=\Big\{ \sup_{0\leq t\leq L^2/2} (\bar Y^y_j(t)-\bar Y^y_j(0))_i \in [L/2, L]\ \ \ \mbox{and}\ \inf_{0\leq t\leq L^2/2} (\bar Y^y_j(t)-\bar Y^y_j(0))_i \in [-L, -L/2] \quad \forall\ 1\leq i\leq d\Big\},
$$
because $|y-z|_\infty \leq L$ by our choice of pairing of $\bar Y^y_j$ and $Y^z_k$. Then by our coupling of $\bar\xi$ and $\xi$, on the event $G_L$, we have
\beq\label{xizetacomp}
\xi(L^2/2, x) \ \geq\ \zeta(x):=\!\!\!\!\!\! \sum_{y\in\Z^d,\atop 1\leq j\leq \bar\xi(0,y)} 1_{E^y_j} 1_{\{\bar Y^y_j(L^2/2)=x\}} \qquad \mbox{for all } |x|_\infty \leq 2L^a.
\eeq
Now observe that, because $(\bar\xi(0,x))_{x\in\Z^d}$ are i.i.d.\ Poisson with mean $\bar\nu$, and $(\bar Y^y_j)_{y\in\Z^d, 1\leq j\leq \bar\xi(0,y)}$ are independent,
$(\zeta(x))_{x\in\Z^d}$ are also i.i.d.\ Poisson distributed with mean $\alpha:=\bar\nu\P^{\bar\xi}(E^y_j)=\bar\nu\P^{\bar\xi}(E^0_1)$, which is bounded away from $0$ uniformly in $L$ by the properties of simple symmetric random walks. This achieves the desired stochastic domination of $\xi$ at time $L^2/2$. Let $\zeta_L(t,\cdot)$ denote the counting field of independent random walks as in (\ref{xi}) with initial condition $\zeta_L(0,y)=\zeta(y)1_{\{|y|_\infty\leq 2L^a\}}$. Then using (\ref{xizetacomp}), uniformly in $x\in\Z^d$ with $|x|_\infty\leq L$, we have
\begin{eqnarray}
&&
\E^\xi\E^X_x\Big[\exp\Big\{-\gamma\int_0^{L^2} \xi(s, X(s))\,{\rm d}s\Big\}\Big]
= \E^{\xi,\bar\xi}\E^X_x\Big[\exp\Big\{-\gamma\int_0^{L^2} \xi(s, X(s))\,{\rm d}s\Big\}\Big] \nn \\
&\leq& \P^{\bar \xi}(G_L^c) + \P^X_x(|X(L^2/2)|_\infty>L^2) + \sup_{|x|_\infty\leq L^2}\E^{\zeta_L}\E^X_x\Big[\exp\Big\{-\gamma\int_0^{L^2/2} \zeta_L(s, X(s))\,{\rm d}s\Big\} \Big] \nn \\
&\leq& 6^d L^{ad} e^{-C_{\nu,\bar\nu} L^d} + e^{-CL^2}+ \sup_{|x|_\infty\leq L^2}\E^{\zeta_L}\E^X_x\Big[\exp\Big\{-\gamma\int_0^{L^2/2} \zeta_L(s, X(s))\,{\rm d}s\Big\} \Big]. \label{pann2}
\end{eqnarray}
By the same argument as for (\ref{urep2}), we have
\beq\label{pann3}
\E^{\zeta_L}\E^X_x\Big[\exp\Big\{-\gamma\int_0^{L^2/2} \zeta_L(s, X(s))\,{\rm d}s\Big\} \Big] = \E^X_x\Big[ \exp \Big\{-\alpha\!\! \sum_{|y|_\infty\leq 2L^a} (1-v_X(L^2/2,y))\Big\}\Big],
\eeq
where
$$
v_X(L^2/2,y)=\E^Y_y\Big[\exp\Big\{-\gamma \int_0^{L^2/2} \delta_0(Y(s)-X(s))\,{\rm d}s\Big\}\Big].
$$
To bound (\ref{pann3}), note that by a union bound in combination with Azuma's inequality we obtain,
\beq \label{pann4}
\sup_{|x|_\infty\leq L^2}\P^X_x\big(\sup_{0\leq s\leq L^2/2}|X(s)|_\infty > 2L^2\big) \leq e^{-CL^2}.
\eeq
On the complementary event $\{\sup_{0\leq s\leq L^2/2}|X(s)|_\infty \leq 2L^2\}$, we have
$$
1-v_X(L^2/2,y) \leq \P^Y_y(\tau_{2L^2} \leq L^2/2) \leq \P({\cal P}_{L^2/2}\geq |y|_\infty-2L^2),
$$
where $\tau_{2L^2}:=\inf\{s\geq 0: |Y(s)|_\infty\leq 2L^2\}$, and ${\cal P}_{L^2/2}$ is a Poisson random variable with mean
$\rho L^2/2$, which counts the number of jumps of $Y$ before time $L^2/2$. Therefore for $L$ sufficiently large,
\begin{eqnarray}
&& \sum_{|y|_\infty> 2L^a} (1-v_X(L^2/2,y)) \leq \sum_{|y|_\infty> 2L^a} \P({\cal P}_{L^2/2}\geq |y|_\infty-2L^2) \nn\\
&\leq&   C \sum_{r=2L^a}^\infty \P({\cal P}_{L^2/2} \geq r/2)\, r^{d-1} \leq C \E[{\cal P}_{L^2/2}^k]\sum_{r=2L^a}^\infty  r^{d-k-1} \nn\\
&\leq& C (\rho L^2/2)^k (2L^a)^{d-k} \leq C L^{-(a-2)k+ad} \leq 1, \label{pann5}
\end{eqnarray}
where we have changed the values of the constant $C$ (independent of $L$) from line to line, and the last inequality holds for all $L$ large if we choose $k$
large enough. Substituting the bounds (\ref{pann4})--(\ref{pann5}) into (\ref{pann3}) then gives the following bound uniformly for $x\in\Z^d$ with $|x|_\infty\leq L^2$:
\begin{eqnarray*}
&& \E^{\zeta_L}\E^X_x\Big[\exp\Big\{-\gamma\int_0^{L^2/2} \zeta_L(s, X(s))\,{\rm d}s\Big\} \Big] \\
&\leq& \P^X_x\big(\sup_{0\leq s\leq L^2/2}|X(s)|_\infty > 2L^2\big) + \E^X_x\Big[ \exp \Big\{-\alpha\!\! \sum_{y\in\Z^d} (1-v_X(L^2/2,y))+1\Big\}\Big] \\
&\leq& e^{-CL^2} + e \E^X_x\Big[ \exp \Big\{-\alpha\!\! \sum_{y\in\Z^d} (1-v_X(L^2/2,y))\Big\}\Big],
\end{eqnarray*}
where by the representation (\ref{urep2}), the expectation is precisely the annealed survival probability of a random walk among a Poisson field of traps with density $\alpha$, for which the bounds in (\ref{annealed}) apply with $\nu$ replaced by $\alpha$ and $t$ by $L^2/2$.
Substituting this bound back into (\ref{pann2}) then gives (\ref{pann}). The same proof applies when we reverse the time direction of $X$ in (\ref{pann}).
\qed

\section{Existence and positivity of the quenched Lyapunov exponent}\label{S:qexistence}
In this section, we prove Theorems~\ref{T:quenched} and \ref{T:quenchedPAM}. In Section~\ref{S:shape}, we state a shape theorem which implies the existence of the quenched Lyapunov exponent for the quenched survival probability $Z^\gamma_{t,\xi}$. In Section~\ref{S:shapepf}, we prove the stated shape theorem for bounded ergodic random fields. In Section~\ref{S:pamexist}, we show how to deduce the existence of the quenched Lyapunov exponent for the solution of the parabolic Anderson model from what we already know for the quenched survival probability. Lastly in Section~\ref{S:quenchedpos}, we prove the positivity of the quenched Lyapunov exponent, which concludes the proof of Theorems~\ref{T:quenched} and \ref{T:quenchedPAM}.

\subsection{Shape theorem and the quenched Lyapunov exponent}\label{S:shape}
In this section, we focus exclusively on the quenched survival probability $Z^\gamma_{t,\xi}$. The approach we adopt in proving the existence of the quenched Lyapunov exponent for $Z^\gamma_{t,\xi}$ uses the subadditive ergodic theorem and follows ideas used by Varadhan in \cite{V03} to prove
the quenched large deviation principle for random walks in random environments.

For $s\geq 0$ and $x\in\Z^d$, let $\P^X_{x,s}$ and $\E^X_{x,s}$ denote respectively probability and expectation for a jump rate $\kappa$ simple symmetric random
walk $X$, starting from $x$ at time $s$. For each $0\le s< t$ and $x,y\in\Z^d$, define
\beq\label{passage}
\begin{aligned}
e(s,t,x,y, \xi) &:=\  \E^X_{x,s}\left[\exp\left\{-\gamma \int_s^t\xi(u,X(u)) \, {\rm d}u\right\}
1_{\{X(t)=y\}}\right], \\
a(s,t,x,y,\xi) &:= \ -\log e(s,t,x,y,\xi).
\end{aligned}
\eeq
We call $a(s,t,x,y,\xi)$ the {\it point to point} passage function from $x$ to $y$ between times $s$ and $t$. We will prove the following shape theorem
for $a(0,t,0,y,\xi)$.

\bt{\bf [Shape theorem]}\label{T:shape}
There exists a deterministic convex function $\alpha: \R^d \to \R$, which we call the shape function, such that $\P^\xi$-a.s., for any compact $K\subset \R^d$,
\beq\label{shape1}
\lim_{t\to\infty} \sup_{y\in tK \cap \Z^d} |t^{-1} a(0,t,0,y,\xi) - \alpha(y/t)| = 0.
\eeq
Furthermore, for any $M>0$, we can find a compact $K\subset \R^d$ such that
\beq\label{shapetight}
\limsup_{t\to\infty} \frac{1}{t} \log \E^X_0\left[\exp\left\{-\gamma \int_0^t\xi(s,X(s))\, {\rm d}s\right\}
1_{\{X(t)\notin tK\}}\right] \leq -M.
\eeq
\et
{\bf Remark.} Note that (\ref{quenched}) in Theorem~\ref{T:quenched} follows easily from Theorem~\ref{T:shape}, which we leave to the reader as an exercise.
In particular, the quenched Lyapunov exponent satisfies
\beq
\tilde\lambda_{d,\gamma,\kappa,\rho, \nu} = \inf_{y\in\R^d} \alpha(y) =\alpha(0) = \lim_{t\to\infty} t^{-1} a(0,t,0,0,\xi),
\eeq
where $\inf_{y\in\R^d} \alpha(y) =\alpha(0)$ follows from the convexity and symmetry  of $\alpha$, since $\xi$ is symmetric.
\medskip

The unboundedness of the random field $\xi$ creates complications for the proof of Theorem~\ref{T:shape}. Therefore we first replace $\xi$
by $\xi_N:=(\max\{\xi(s,x), N\})_{s\geq 0, x\in\Z^d}$ for some large $N>0$ and prove a corresponding shape theorem, then use almost sure
properties of $\xi$ established by Kesten and Sidoravicius in~\cite{KS03} to control the error caused by the truncation.

\bt{\bf [Shape theorem for bounded ergodic potentials]}\label{T:shapeb}
Let $\zeta:=(\zeta(s,x))_{s\geq 0, x\in\Z^d}$ be a real-valued random field which is ergodic with respect to the shift map
$\theta_{r,z}\zeta := (\zeta(s+r, x+z))_{s\geq 0, x\in\Z^d}$, for all $r\geq 0$ and $z\in\Z^d$. Assume further that $|\zeta(0,0)|\leq A$ a.s.\ for
some $A>0$. Then the conclusions of Theorem~\ref{T:shape} hold with $\xi$ replaced by $\zeta$.
\et
{\bf Remark.} Note that Theorem~\ref{T:shapeb} can be applied to the occupation field of the exclusion process or the voter
model in an ergodic equilibrium, which in particular implies the existence of the corresponding quenched Lyapunov exponents.
\medskip

Before we prove Theorem~\ref{T:shapeb} in the next section, let us first show how to deduce Theorem~\ref{T:shape} from Theorem~\ref{T:shapeb}, using
almost sure bounds on $\xi$ from \cite{KS03}.
\bigskip

\noindent
{\bf Proof of Theorem~\ref{T:shape}.} Note that, since $\xi$ is non-negative, (\ref{shapetight}) follows from elementary large deviation estimates
for the random walk $X$, if we take $K$ to be a large enough closed ball centered at the origin, which we fix for the rest of the proof.

By applying Theorem~\ref{T:shapeb} to the truncated random field $\xi_N$, we have that for each $N>0$, there exists a convex shape function
$\alpha_N :\R^d\to\R$ such that (\ref{shape1}) holds with $\xi$ replaced by $\xi_N$ and $\alpha$ replaced by $\alpha_N$. Note that $\alpha_N$ is
monotonically increasing in $N$, and its limit $\alpha$ is necessarily convex. To prove (\ref{shape1}), it then suffices to show that, for any $\eps>0$,
we can choose $N$ sufficiently large such that $\P^\xi$-a.s.,
\beq\label{shapepf1}
\frac{1}{t} \sup_{y\in tK\cap \Z^d} |a(0,t,0,y,\xi) - a(0,t,0,y,\xi_N)|\leq \eps \qquad \mbox{for all $t$ sufficiently large}.
\eeq

To prove (\ref{shapepf1}), we will need Lemma 15 from \cite{KS03}, which by Borel-Cantelli implies that there exist positive constants $C_0, C_1, C_2, C_3, C_4$ with $C_0>1$,
such that if $\Xi_l$ denotes the space of all possible random walk trajectories $\pi: [0,t]\to \Z^d$, which contain exactly $l$ jumps and are contained
in the rectangle $[-C_1 t\log t, C_1 t\log t]^d$, then $\P^\xi$-a.s., for all $t\in\N$ sufficiently large, we have
\beq\label{shapepf2}
\sup_{\pi \in \Xi_l} \int_0^t \xi(s, \pi(s)) 1_{\{\xi(s,\pi(s))\geq C_2\nu C_0^{dm}\}}\, {\rm d}s \leq
(t+l) \sum_{r=m}^\infty C_3  C_0^{r(d+6)+d} e^{-C_4 C_0^{r/4}} \qquad \forall\ m\in\N, l\geq 0,
\eeq
where $A_m:=\sum_{r=m}^\infty C_3  C_0^{r(d+6)+d} e^{-C_4 C_0^{r/4}}\to 0$ as $m\to\infty$.

One important consequence of (\ref{shapepf2}) is that
\beq\label{shapepf3}
0<\sup_{y\in K} \alpha(y) <\infty.
\eeq
Indeed, if $l_t(X)$ denotes the number of jumps of $X$ on the time interval $[0,t]$, then
\begin{eqnarray*}
\sup_{y\in K} \alpha(y) &\leq& \lim_{t\to\infty} -t^{-1} \log \inf_{y\in tK\cap \Z^d}\E^X_0\left[\exp\left\{-\gamma \int_0^t\xi(s,X(s))\, {\rm d}s\right\}
1_{\{X(t)=y\}}\right] \\
&\leq& \lim_{t\to\infty} -t^{-1} \log \inf_{y\in tK\cap \Z^d} \E^X_0\left[\exp\left\{-\gamma \int_0^t\xi(s,X(s))\, {\rm d}s\right\}
1_{\{X(t)=y, l_t(X)\leq 2D(K)t, X\in \Xi_{l_t(X)}\}}\right],
\end{eqnarray*}
where $D(K):=\sup_{y\in K} \vert y\vert_1$. We can then apply (\ref{shapepf2}) and large deviation estimates for random walks to the above bound to deduce
$\sup_{y\in K} \alpha(y) <\infty$. The fact that $\sup_{y\in K} \alpha(y)>0$ for a large ball $K$ again follows from basic large deviation estimates.

By large deviation estimates, we can find $B$ large enough such that
\beq\label{shapepf4}
\P^X_0(l_t(X)\geq Bt\ \mbox{or} \ X\notin \Xi_{l_t(X)})\leq e^{-2 \sup_{y\in K} \alpha(y)\, t} \qquad \mbox{for all $t$ sufficiently large}.
\eeq
Let $N=C_2\nu C_0^{dm}$. Then by (\ref{shapepf2}), $\P^\xi$-a.s., uniformly in $y\in\Z^d$ and for all $t$ large, we have
\begin{eqnarray}
e(0,t,0,y, \xi) &\geq& e^{-(1+B)A_m\gamma t} \E^X_0\left[\exp\left\{-\gamma \int_0^t\xi_N(s,X(s))\, {\rm d}s\right\} 1_{\{X(t)=y, l_t(X)\leq Bt, X\in \Xi_{l_t(X)}\}}\right] \nonumber \\
&\geq& e^{-(1+B)A_m\gamma t} (e(0,t,0,y,\xi_N)- e^{-2 \sup_{y\in K} \alpha(y)\, t}), \label{shapepf5}
\end{eqnarray}
where in the last inequality we applied (\ref{shapepf4}). Since $-t^{-1}\log e(t,0,y,\xi_N)\to\alpha_N(y/t)$ uniformly for $y\in tK\cap \Z^d$ by Theorem~\ref{T:shapeb}, and $\sup_{y\in K}\alpha_N(y)\leq \sup_{y\in K}\alpha(y)$, (\ref{shapepf5}) implies that $\P^\xi$-a.s., uniformly in $y\in tK\cap \Z^d$ and for all $t$ large, we have
$$
t^{-1}a(0,t,0,y, \xi) \leq t^{-1} a(0,t,0,y, \xi_N) + (1+B)A_m\gamma +o(1).
$$
Since $a(0,t,0,y,\xi) \geq a(0,t,0,y,\xi_N)$, and $A_m$ can be made arbitrarily small by choosing $m$ sufficiently large, (\ref{shapepf1}) then follows.
\qed

\bigskip
\noindent
{\bf Remark.} Theorem~\ref{T:shape} in fact holds for the catalytic case as well, where we take $\gamma<0$ in (\ref{passage}) and (\ref{shapetight}). This
implies the existence of the quenched Lyapunov exponent in Theorem~\ref{T:quenched} for the catalytic case, where we set $\gamma<0$ in the definition of $Z^\gamma_{t,\xi}$. Indeed, Theorem~\ref{T:shapeb} still applies to the truncated field $\xi_N$. To control the error caused by the truncation, the following
modifications are needed in the proof of Theorem~\ref{T:shape}. To prove (\ref{shapetight}), we need to apply (\ref{shapepf2}). More precisely, we
need to first consider trajectories $(X_x)_{0\leq s\leq t}$ which are not contained in $[-C_1 t\log t, C_1 t\log t]^d$. The contribution from these trajectories
can be shown to decay super-exponentially in $t$ by large deviation estimates and a bound on $\xi$ given in (2.37) of \cite[Lemma 4]{KS03}. For $X$ which lies inside $[-C_1 t\log t, C_1 t\log t]^d$, we can then use (\ref{shapepf2}) and large deviations to deduce (\ref{shapetight}). In contrast to (\ref{shapepf5}), we need to upper bound $e(0,t,0,y,\xi)$ in terms of $e(0,t,0,y,\xi_N)$. The proof is essentially the same, except that in place of (\ref{shapepf4}), we need to show that we can choose
$B$ large enough, such that $\P^\xi$-a.s.,
\beq
\sup_{y\in tK\cap\Z^d}\!\!\E^X_0\left[\exp\left\{|\gamma| \int_0^t\xi(s,X(s))\, {\rm d}s\right\} 1_{\{X(t)=y\}} 1_{\{l_t(X)\geq Bt \mbox{ or } X\notin \Xi_{l_t(X)}\}}\right]
\leq \!\! \inf_{y\in tK\cap\Z^d}\!\! \P^X_0(X_t=y).
\eeq
This can be proved by appealing to (\ref{shapetight}), and applying (\ref{shapepf2}) and large deviation estimates.

\subsection{Proof of shape theorem for bounded ergodic potentials}\label{S:shapepf}

In this section, we prove Theorem~\ref{T:shapeb}. From now on, let $\Q_+$ denote the set of positive rationals, and let
$\Q^d$ denote the set of points in $\R^d$ with rational coordinates. We start with the following auxiliary result.

\blem\label{L:subadd}  There exists a deterministic function $\alpha: \Q^d\to [-\gamma A,\infty)$ such that for every $y\in \Q^d$,
\beq\label{subadd1}
\lim_{t\to\infty \atop ty\in\Z^d} t^{-1} a(0,t,0,ty, \zeta)=\alpha(y) \qquad \P^\zeta-a.s.
\eeq
\elem
{\bf Proof.} Since we assume $y\in\Q^d$ and $ty\in\Z^d$ in (\ref{subadd1}), without loss of generality, it suffices to consider $y\in\Z^d$ and $t\in\N$.
Note that by the definition of the passage function $a$ in (\ref{passage}), $\P^\zeta$-a.s.,
\beq\label{triangleineq}
a(t_1,t_3,x_1,x_3, \zeta) \leq a(t_1, t_2, x_1, x_2, \zeta) + a(t_2, t_3, x_2, x_3, \zeta) \quad \forall\, t_1<t_2<t_3, \, x_1, x_2, x_3\in\Z^d.
\eeq
Together with our assumption on the ergodicity of $\zeta$, this implies that the two-parameter family $a(s,t,sy,ty,\zeta)$, $0\leq s\leq t$ with $s,t\in\Z$,
satisfies the conditions of Kingman's subadditive ergodic theorem (see e.g.~\cite{L85}). Therefore, there exists a deterministic constant $\alpha(y)$ such that (\ref{subadd1}) holds. The fact that $\alpha(y)\geq -\gamma A$ follows by bounding $\zeta$ from above by the uniform bound $A$, and $\alpha(y)<\infty$ follows
from large deviation bounds for the random walk $X$.
%
\qed
\medskip

To extend the definition of $\alpha(y)$ in Lemma~\ref{L:subadd} to $y\notin \Q^d$ and to prove the uniform convergence in (\ref{shape1}), we need to establish
equicontinuity of $t^{-1}a(0,t,0,ty, \zeta)$ in $y$, as $t\to\infty$. For that, we first need a large deviation estimate for the random walk $X$.

\blem\label{L:rwldp}
Let $X$ be a jump rate $\kappa$ simple symmetric random walk on $\Z^d$ with $X(0)=0$. Then for every $t>0$ and $x\in\mathbb Z^d$, we have
\beq\label{rwldp}
\P_0^X(X(t)=x)= \frac{e^{-J(\frac{x}{t})\,t}}{(2\pi t)^{\frac{d}{2}}\Pi_{i=1}^d
\big(\frac{x_i^2}{t^2}+\frac{\kappa^2}{d^2}\big)^{1/4}} \left(1+o(1)\right),
\eeq
where
$$
J(x) := \sum_{i=1}^d \frac{\kappa}{d} j\Big(\frac{d x_i}{\kappa}\Big) \qquad \mbox{with} \quad  j(y) := y\sinh^{-1} y -\sqrt{y^2+1}+1,
$$
and the error term $o(1)$ tends to zero as $t\to\infty$ uniformly in $x\in tK\cap \Z^d$, for any compact $K\subset \R^d$.
\elem
{\bf Proof.} Since the coordinates of $X$ are independent, it suffices to consider the case $X$ is a rate $\kappa/d$ simple symmetric random
walk on $\Z$. Let $\sigma:=t/\lceil t\rceil$. Let $Z^\lambda_1,\cdots, Z^\lambda_{\lceil t\rceil}$ be i.i.d.\
with
$$
\P(Z^\lambda_1 = y) = \P(X(\sigma)=y) e^{\lambda y -\Phi(\lambda)}, \qquad y\in\Z,
$$
where
$$
\Phi(\lambda) = \log \E[e^{\lambda X(\sigma)}] = \frac{\sigma \kappa}{d}(\cosh \lambda -1).
$$
Note that
$$
\E[Z^\lambda_1] =  \frac{\rm d\Phi}{{\rm d}\lambda}(\lambda) = \frac{\sigma \kappa}{d}\sinh \lambda
\qquad \mbox{and} \qquad \mbox{Var}(Z^\lambda_1) = \frac{\rm d^2\Phi}{{\rm d}^2\lambda}(\lambda) = \frac{\sigma \kappa}{d}\cosh \lambda.
$$
We shall set $\lambda=\sinh^{-1}(\frac{d x}{\kappa t})$ so that $\E[Z^\lambda_1]=x/\lceil t\rceil$. If we let $S_{\lceil t\rceil}:=\sum_{i=1}^{\lceil t\rceil} Z^\lambda_i$, then observe that
$$
\P_0^X(X(t)=x)=\P(S_{\lceil t\rceil}=x)e^{-\lambda x+\lceil t\rceil\Phi(\lambda)} = \P(S_{\lceil t\rceil}=x)e^{-\frac{\kappa}{d} j(\frac{d x}{\kappa t}) t}.
$$
Note that $S_{\lceil t\rceil}-x$ has mean $0$, variance $t\sqrt{\frac{x^2}{t^2}+\frac{\kappa^2}{d^2}}$, and characteristic function
$$
e^{\lceil t\rceil\big(\Phi(ik+\lambda)-\Phi(\lambda)\big)-ikx}=e^{ix(\sin k-k)-t\sqrt{\frac{x^2}{t^2}+\frac{\kappa^2}{d^2}}(1-\cos k)}.
$$
Applying Fourier inversion then gives (\ref{rwldp}).
\qed

With the help of Lemma \ref{L:rwldp}, we can control the modulus of continuity of $t^{-1}a(0,t,0,ty, \zeta)$.
\blem\label{L:modcont}
Let $K$ be any compact subset of $\R^d$. There exists $\phi_K: (0,\infty)\to(0,\infty)$ with $\lim_{r\downarrow 0}\phi_K(r)=0$, such that
for any $\eps>0$, $\P^\zeta$-a.s., we have
\beq\label{modcont1}
\limsup_{t\to\infty} \sup_{x,y\in tK \cap \Z^d \atop \Vert x-y\Vert \leq \eps t}
t^{-1}|a(0,t,0,x,\zeta)-a(0,t,0,y,\zeta)| \leq \phi_K(\eps).
\eeq
\elem
{\bf Proof.} Let $K\subset \R^d$ be compact. It suffices to consider $\eps\in (0, 1/2)$, which we also fix from now on.
First note that, by Lemma \ref{L:rwldp}, $\P^\zeta$-a.s.,
\beq\label{modcont2}
\inf_{z\in tK\cap \Z^d} e(0,t,0,z,\zeta) \geq e^{-At} \inf_{z\in tK\cap \Z^d}\P^X_0(X_t=z) \geq e^{-(A+1)t-\sup_{u\in K} J(u)\, t}
\eeq
for all $t$ sufficiently large.

Also note that for all $z\in\Z^d$ and $t>0$,
\beq\label{modcont3}
e(0,t,0,z, \zeta) = \sum_{w\in\Z^d} e(0,(1-\eps)t, 0, w,\zeta)\, e((1-\eps)t, t, w, z, \zeta).
\eeq
By large deviation estimates, we can choose a ball $B_R$ centered at the origin with radius $R$ large enough and independent of $\eps$,
such that $K\subset B_R$ and $\P^\zeta$-a.s.,
\begin{eqnarray}
&& \sup_{z\in\Z^d} \sum_{w\in \Z^d \atop w\notin tB_R} e(0,(1-\eps)t, 0, w,\zeta)\, e((1-\eps)t, t, w, z, \zeta) \nonumber \\
&\leq& \P^X_0\big(X((1-\eps)t)\notin tB_R\big) e^{At} \leq e^{-(A+2)t-\sup_{u\in K} J(u)\, t} \nonumber
\end{eqnarray}
for all $t$ sufficiently large. In view of (\ref{modcont2}), the dominant contribution in (\ref{modcont3}) comes from $w\in tB_R\cap \Z^d$.
Therefore to prove (\ref{modcont1}), it suffices to verify
\beq\label{modcont4}
\limsup_{t\to\infty} \sup_{x,y\in tK \cap \Z^d \atop \Vert x-y\Vert \leq \eps t}
t^{-1}\left|\log \frac{\sum_{w\in tB_R\cap\Z^d} e(0,(1-\eps)t, 0, w,\zeta)\, e((1-\eps)t, t, w, y, \zeta)}{\sum_{w\in tB_R\cap\Z^d} e(0,(1-\eps)t, 0, w,\zeta)\, e((1-\eps)t, t, w, x, \zeta)}\right| \leq \phi_K(\eps).
\eeq
Note that $\P^\zeta$-a.s., and uniformly in $x,y\in tK\cap \Z^d$ with $\Vert x-y\Vert \leq \eps t$,
\begin{eqnarray*}
&& \frac{\sum_{w\in tB_R\cap\Z^d} e(0,(1-\eps)t, 0, w,\zeta)\, e((1-\eps)t, t, w, y, \zeta)}{\sum_{w\in tB_R\cap\Z^d} e(0,(1-\eps)t, 0, w,\zeta)\, e((1-\eps)t, t, w, x, \zeta)} \\
&\leq& \sup_{w\in tB_R\cap \Z^d} \frac{e((1-\eps)t, t, w, y, \zeta)}{e((1-\eps)t, t, w, x, \zeta)}
\leq \sup_{w\in tB_R\cap \Z^d} \frac{\P^X_0(X(\eps t)=y-w)}{\P^X_0(X(\eps t)=x-w)} \ e^{2A\eps t} \\
&\leq& \exp\left\{\eps t \sup_{w\in tB_R\cap \Z^d}\Big(J\Big(\frac{x-w}{\eps t}\Big) -J\Big(\frac{y-w}{\eps t}\Big)\Big)+3A\eps t\right\} \\
&\leq& \exp\left\{\eps t \sup_{u,v\in B_{2R/\eps}, \Vert u-v\Vert\leq 1} |J(u)-J(v)| +3A\eps t\right\}
\end{eqnarray*}
for all $t$ sufficiently large, where we applied Lemma~\ref{L:rwldp}, and $B_{2R/\eps}$ denotes the ball of radius $2R/\eps$, centered at the origin.
Therefore (\ref{modcont4}) holds with
$$
\phi_K(\eps) = 3A\eps + \eps \sup_{u,v\in B_{2R/\eps}, \Vert u-v\Vert\leq 1} |J(u)-J(v)|.
$$
It only remains to verify that $\phi_K(\eps)\downarrow 0$ as $\eps\downarrow 0$, which is easy to check from the definition of $J$.
\qed
\bigskip

\noindent
{\bf Proof of Theorem~\ref{T:shapeb}.} Because $\zeta$ is uniformly bounded, (\ref{shapetight}) follows by large deviation estimates for the number of
jumps of $X$ up to time $t$. Lemma \ref{L:modcont} implies that for each compact $K\subset\R^d$, the function $\alpha$ in Lemma~\ref{L:subadd} satisfies
\beq\label{shapebpf0}
\sup_{u,v\in K\cap \Q^d\atop \Vert u-v\Vert \leq \eps} |\alpha(u)-\alpha(v)|\leq \phi_K(\eps) \qquad \mbox{for all }\ \eps>0.
\eeq
This allows us to extend $\alpha$ to a continuous function on $\R^d$.

To prove (\ref{shape1}), it suffices to show that for each $\delta>0$,
\beq\label{shapebpf1}
\limsup_{t\to\infty} \sup_{y\in tK \cap \Z^d} |t^{-1} a(0,t,0,y,\xi) - \alpha(y/t)| \leq \delta.
\eeq
We can choose an $\eps$ such that $\phi_K(\eps)<\delta /3$. We can then find a finite number of points $x_1,\cdots, x_m\in \Q^d$ which form an $\eps$-net in $K$,
and along a subsequence of times of the form $t_n=n\sigma$ with $\sigma x_i\in\Z^d$ for all $x_i$, we have $t_n^{-1}a(0,t_n,0,t_n x_i)\to\alpha(x_i)$ a.s.
The uniform control of modulus of continuity provided by Lemma~\ref{L:modcont} and (\ref{shapebpf0}) then implies (\ref{shapebpf1}) along $t_n$. This can be transferred
to $t\to\infty$ along $\R$ using
$$
e(0,t,0,y,\zeta) \geq e(0,s,0,y,\zeta)\, e(s,t,y,y,\zeta) \geq e(0,s,0,y,\zeta) e^{-(\kappa+\gamma A)(t-s)} \quad \mbox{for } s<t.
$$
Lastly, to prove the convexity of $\alpha$, let $x,y\in \R^d$ and $\beta\in (0,1)$. Then $\P^\zeta$-a.s., we have
$$
a(0,t_n,0, \beta y_n+(1-\beta) x_n, \zeta) \leq a(0, \beta t_n, 0, \beta y_n, \zeta) + a(\beta t_n, t_n, \beta y_n, \beta y_n+(1-\beta)x_n, \zeta),
$$
where we take sequences $t_n,x_n, y_n$ with $t_n\to\infty$, $x_n/t_n\to x$, $y_n/t_n\to y$, and $\beta y_n, (1-\beta)x_n\in\Z^d$. By Lemma~\ref{L:subadd},
the first term divided by $t_n$ converges a.s.\ to $\alpha(\beta y+(1-\beta)x)$, the second term divided by $\beta t_n$
converges a.s.\ to $\alpha(y)$, while the last term divided by $(1-\beta)t_n$ converges in probability to $\alpha(x)$ by translation invariance. The convexity
of $\alpha$ then follows.
\qed

\subsection{Existence of the quenched Lyapunov exponent for the PAM}\label{S:pamexist}

{\bf Proof of (\ref{quenchedPAM}) in Theorem~\ref{T:quenchedPAM}.} Since $Z^\gamma_{t,\xi}$ is equally distributed with $u(t,0)$
for each $t\geq 0$, $-t^{-1}\log u(t,0)$ converges in probability to the quenched Lyapunov exponent $\tilde\lambda_{d,\gamma,\kappa,\rho, \nu}$. It only remains to verify the almost sure convergence. We will bound the variance of $\log u(t,0)$, which is the same as that of $\log Z^\gamma_{t,\xi}$, and then apply Borel-Cantelli.

Assume that $t\in\N$. Note that we can write $\xi$ as a sum of i.i.d.\ random fields $(\xi_i(s,x))_{s\geq 0, x\in\Z^d}$, each of which is defined
from a Poisson system of independent random walks with density $\nu/t$, in the same way as $\xi$. Then we can perform a martingale decomposition and write
$$
\log Z^\gamma_{t,\xi} - \E^\xi[\log Z^\gamma_{t,\xi}] = \sum_{i=1}^t V_i :=\sum_{i=1}^t \Big(\E^\xi[\log Z^\gamma_{t,\xi} |\xi_1,\cdots, \xi_i] - \E^\xi[\log Z^\gamma_{t,\xi} |\xi_1,\cdots, \xi_{i-1}]\Big),
$$
and hence ${\rm Var}(\log Z^\gamma_{t,\xi}) = \sum_{i=1}^t \E^\xi[V_i^2]$.

For each $1\leq i\leq t$, we have
\begin{eqnarray*}
V_i &=& \E^{\xi_{i+1},\cdots, \xi_t}\big[\log Z^\gamma_{t,\xi} - \E^{\xi_i}[\log Z^\gamma_{t,\xi}] \big] \\
&=& \E^{\xi_{i+1},\cdots, \xi_t} \E^{\xi_i'}\left[\log \frac{\E^X_0\big[e^{-\gamma \int_0^t \big(\sum_{1\leq j\leq t, j\neq i} \xi_j(s, X(s)) + \xi_i(s, X(s))\big)\, {\rm d}s}\big]}{\E^X_0\big[e^{-\gamma \int_0^t \big(\sum_{1\leq j\leq t, j\neq i} \xi_j(s, X(s)) + \xi'_i(s, X(s))\big)\, {\rm d}s}\big]} \right] \\
&=& \E^{\xi_{i+1},\cdots, \xi_t} \E^{\xi_i'}\left[\log \E^{X,i}\big[e^{-\gamma \int_0^t \xi_i(s,X(s))\, {\rm d}s}\big] -\log \E^{X,i}\big[e^{-\gamma \int_0^t \xi'_i(s,X(s))\, {\rm d}s} \big] \right],
\end{eqnarray*}
where $\xi_i'$ denotes an independent copy of $\xi_i$, and $\E^{X, i}$ denotes expectation with respect to the Gibbs transform of the
random walk path measure $\P^X_0$, with Gibbs weight $e^{-\gamma \int_0^t \sum_{1\leq j\leq t, j\neq i} \xi_j(s, X(s))\, {\rm d}s}$. Then by Jensen's inequality,
\begin{eqnarray}
\E^\xi [V_i^2] &\leq& \E^{\xi, \xi_i'} \Big[\big(\log \E^{X,i}\big[e^{-\gamma \int_0^t \xi_i(s,X(s))\, {\rm d}s}\big] -\log \E^{X,i}\big[e^{-\gamma \int_0^t \xi'_i(s,X(s))\, {\rm d}s} \big]\big)^2 \Big] \nonumber\\
&\leq& 2 \E^{\xi, \xi_i'} \Big[\big(\log \E^{X,i}\big[e^{-\gamma \int_0^t \xi_i(s,X(s))\, {\rm d}s}\big]\big)^2\Big] + 2\E^{\xi, \xi_i'} \Big[\big(\log \E^{X,i}\big[e^{-\gamma \int_0^t \xi'_i(s,X(s))\, {\rm d}s} \big]\big)^2 \Big] \nonumber \\
&=& 4 \E^\xi \Big[\big(\log \E^{X,i}\big[e^{-\gamma \int_0^t \xi_i(s,X(s))\, {\rm d}s}\big]\big)^2\Big] \nonumber \\
&\leq& 4 \E^\xi \Big[\Big(\E^{X,i}\Big[\gamma \int_0^t \xi_i(s,X(s))\, {\rm d}s\Big]\Big)^2\Big] \nonumber\\
&\leq& 4 \gamma^2 \E^\xi \E^{X,i}\Big[\Big(\int_0^t \xi_i(s,X(s))\, {\rm d}s\Big)^2\Big] = 4 \gamma^2 \E^{X,i}\E^{\xi_i}\Big[\Big(\int_0^t \xi_i(s,X(s))\, {\rm d}s\Big)^2\Big] , \nonumber 
\end{eqnarray}
where in the third line we used the exchangeability of $\{\xi_i, \xi_i'\}$, and in the fourth line we applied Jensen's inequality\footnote{Note that this is where the proof fails for the $\gamma<0$ case.} to the non-negative convex function $-\log x$ on the interval $(0,1]$.

Note that for any realization of $((X(s))_{0\leq s\leq t}$, we have
\begin{eqnarray*}
&& \E^{\xi_i}\Big[\Big(\int_0^t \xi_i(s,X(s))\, {\rm d}s\Big)^2\Big] = 2\iint\limits_{0<u<v<t} \E^{\xi_i}[\xi_i(u, X(u))\xi_i(v,X(v))]\, {\rm d}u\, {\rm d}v \\
&=& 2 \iint\limits_{0<u<v<t} \Big(\frac{\nu^2}{t^2}\sum_{y\in\Z^d \atop y\neq X(u)} \P^Y_{y,u}(Y(v)=X(v)) +  \big(\frac{\nu^2}{t^2}+\frac{\nu}{t}\big)\P^Y_{X(u),u}(Y(v)=X(v))\Big) \, {\rm d}u\, {\rm d}v \\
&\leq& 2\nu^2 + 2\nu \int_0^t \P^Y_{0,0}(Y(s)=0) \,{\rm d}s,
\end{eqnarray*}
where $\P^Y_{y,s}$ denotes probability for a simple symmetric random walk on $\Z^d$ with jump rate $\rho$, starting from $y$ at time $s$, and
in the last line we used that $\P^Y_{0,0}(Y(s)=y)$ is maximized at $y=0$ for all $s\geq 0$. Combined with the previous bounds, we obtain
$$
{\rm Var}(\log u(t,0)) = {\rm Var}(\log Z^\gamma_{t,\xi}) = \sum_{i=1}^t \E^\xi[V_i^2] \leq 8\gamma^2\nu^2 t + 8\gamma^2 \nu t \int_0^t \P^Y_{0,0}(Y(s)=0) \, {\rm d}s \leq Ct^{\frac{3}{2}}
$$
for some $C>0$, since $\int_0^t \P^Y_{0,0}(Y(s)=0) \, {\rm d}s$ is of order $\sqrt{t}$ in dimension $d=1$, of order $\log t$ in $d=2$, and converges in $d\geq 3$. Therefore for any $\eps>0$,
$$
\P^\xi\big( \big|\log u(t,0) - \E^\xi[\log u(t,0)] \big| \geq \eps t\big) \leq \frac{C}{\eps^2 \sqrt{t}},
$$
which by Borel-Cantelli implies that along the sequence $t_n=n^3$, $n\in\N$, we have almost sure convergence of $-t^{-1}\log u(t,0)$ to the quenched Lyapunov exponent
$\tilde\lambda_{d,\gamma,\kappa,\rho, \nu}$.

To extend the almost sure convergence to $t\to\infty$ along $\R$, consider $t\in [t_n, t_{n+1})$ for some $n\in\N$. As at the end of the proof of Proposition
\ref{P:subexp}, we have
$$
\begin{aligned}
u(t,0) &\geq  e^{-\kappa(t-t_n)} e^{-\gamma \int_{t_n}^t \xi(s,0)\, {\rm d}s} u(t_n,0),  \\
u(t,0) &\leq e^{\kappa(t_{n+1}-t)} e^{\gamma \int_t^{t_{n+1}} \xi(s,0)\, {\rm d}s} u(t_{n+1},0).
\end{aligned}
$$
Note that $(t_{n+1}-t_n)/t_n\to 0$ as $n\to\infty$, and we claim that also $t_n^{-1}\int_{t_n}^{t_{n+1}} \xi(s,0)\, {\rm d}s \to 0$ a.s.\ as $n\to\infty$, which
then implies the desired almost sure convergence of $t^{-1}\log u(t,0)$ as $t\to\infty$ along $\R$. Indeed, since $\int_0^1 \xi(s,0)\, {\rm d}s$ has finite
exponential moments, as can be seen from (\ref{urep3}) applied to the case $\gamma<0$ and $X\equiv 0$, we have exponential tail bounds on $\int_0^1 \xi(s,0)\, {\rm d}s$, which by Borel-Cantelli implies that a.s.\ $\sup_{0\leq i< m}\int_{i}^{i+1}\xi(s,0) \, {\rm d}s \leq \log m$ for all $m\in\N$ sufficiently
large. The above claim then follows.
\qed

\subsection{Positivity of the quenched Lyapunov exponent}\label{S:quenchedpos}
In this section, we conclude the proof of Theorems \ref{T:quenched} and \ref{T:quenchedPAM} by showing that the quenched Lyapunov exponent $\tilde\lambda_{d,\gamma,\kappa,\rho, \nu}$ is positive in all dimensions. The strategy is as follows: Employing a result of Kesten and Sidoravicius~\cite[Prop.~8]{KS05}, we deduce that $\P^\xi$-a.s. for eventually all integer time points $t$, sufficiently many $X$ paths encounter a $\xi$-particle close-by for of order $t$ many integer time points.  Using the Markov property, we then show that with positive $\P_0^X$ probability, $X$ moves to a close-by $\xi$-particle (which itself stays at its site for some time)
within a very short time interval and collects some local time with this $\xi$-particle. This then implies the desired exponential decay.
\bigskip

\noindent{\bf Proof of Theorems \ref{T:quenched} and \ref{T:quenchedPAM}.} Since we have shown the quenched Lyapunov exponent $\tilde\lambda_{d,\gamma,\kappa,\rho, \nu}$ in Theorems~\ref{T:quenched} and \ref{T:quenchedPAM} to be the same, it suffices to consider only Theorem~\ref{T:quenched}. Note that the upper bound on $\tilde\lambda_{d,\gamma,\kappa,\rho, \nu}$ in Theorem~\ref{T:quenched} follows trivially by requiring the walk $X$ to stay at the origin. To show $\tilde\lambda_{d,\gamma,\kappa,\rho, \nu}>0$, we will make the strategy outlined above precise. In compliance with \cite{KS05} we let $C_0$ and $r > 0$ be large integers
and for $\vec i \in \Z^d$ define the cubes
$$
\mathcal{Q}_r(\vec i) := \prod_{j=1}^d [i_j, i_j + C_0^r).
$$
In a slight abuse of common notation, let $D([0,\infty),\Z^d)$ denote the Skorohod space restricted to those functions that
start in $0$ at time $0$ and have
nearest neighbour jumps only.
Then set
\begin{align*}
J_k := \big\{ \Phi \in D([0,\infty),\Z^d) : \Phi \text{ jumps at most } dC_0^r(\kappa \vee 1)k \text{ times up to time } k \big\}.
\end{align*}
For integer times $t > 0$ define
\begin{align*}
\Xi(t) := \bigcap_{k=\lfloor t/4 \rfloor}^t J_k.
\end{align*}
Then standard large deviation bounds yield
\begin{equation} \label{eq:manyJumpsEst}
\P_0^X \big( X \in \Xi(t)^c \big) \leq e^{-c(t + o(t))},
\end{equation}
for some $c > 0.$
In addition, define the cube
$$
\mathcal{C}_t := [-dC_0^r(\kappa \vee 1)t, dC_0^r(\kappa \vee 1)t]^d \cap \Z^d,
$$
as well as for arbitrary $t \in \N,$ $k \in \{0, \cdots, t\},$ $\Phi \in \Xi(t)$ and $\epsilon \geq 0$ the events
$$
A(t,\Phi,k, \epsilon) := \big \{\exists \vec i \in \mathcal{C}_t \,: \, \Phi(k) \in \mathcal{Q}_r(\vec i) \text{ and }
\exists \, y \in \mathcal{Q}_r(\vec i) \, : \, \xi(s,y) \geq 1 \, \forall \, s \in [k,k+\epsilon/\rho]\big \}
$$
and
\begin{equation*}
G(t) := \bigcap_{\Phi \in \Xi(t)} \Big\{\sum_{k \in \{\lfloor t/4 \rfloor,\cdots, t-1\}}
1_{A(t,\Phi,k, \epsilon)} \geq \epsilon t \Big\},
\end{equation*}
which both depend on $\xi.$

For $\epsilon$ small enough, using Borel-Cantelli, it is a consequence of \cite[Prop.~8]{KS05} that $\P^\xi$-a.s., $G(t)$ occurs for
eventually all $t \in \N.$
Indeed, denoting by
$
\Xi (t)|_{\{\lfloor t/4 \rfloor, \cdots, t\}}
$
the subset of $(\Z^d)^{\{\lfloor t/4 \rfloor, \cdots, t\}}$ obtained by restricting each element of
$\Xi(t)$ to the domain ${\{\lfloor t/4 \rfloor, \cdots, t\}}$, we estimate
\begin{eqnarray}
 \P^\xi \big( G(t)^c \big)
\!\!\!&\leq&\!\!\! \P^\xi \Big( \bigcup_{\Phi \in \Xi(t)}
\Big\{ \sum_{k \in \{\lfloor t/4 \rfloor,\cdots, t-1\}}
1_{A(t,\Phi,k, \epsilon)} \leq \epsilon t \Big\} \Big) \nonumber\\
\!\!\!&\leq&\!\!\! \P^\xi \Big( \bigcup_{\Phi \in \Xi(t)}
\Big\{ \sum_{k \in \{\lfloor t/4 \rfloor,\cdots, t-1\}}
1_{A(t,\Phi,k, 0)} \leq \frac{t}{2} \Big\} \Big) \nonumber \\
\!\!\! && + \big\vert \Xi (t)|_{\{\lfloor t/4 \rfloor, \cdots, t\}} \big\vert
\times \max_{\Phi \in \Xi(t)} \P^\xi \Big(
\sum_{k \in \{\lfloor t/4 \rfloor,\cdots, t-1\}} 1_{A(t,\Phi,k, 0)} \geq \frac{t}{2}, \sum_{k \in \{\lfloor t/4 \rfloor,\cdots, t-1\}}
1_{A(t,\Phi,k, \epsilon)} \leq \epsilon t \Big) \nonumber\\
\!\!\!&\leq&\!\!\!
\P^\xi \Big( \bigcup_{\Phi \in \Xi(t)}
\Big\{ \sum_{k \in \{\lfloor t/4 \rfloor,\cdots, t-1\}}
1_{A(t,\Phi,k, 0)} \leq \frac{t}{2} \Big\} \Big)
+ \big\vert \Xi (t)|_{\{\lfloor t/4 \rfloor, \cdots, t\}} \big\vert \times
\P \Big( \sum_{i=1}^{t/2} p_{i,\epsilon} \leq \epsilon t \Big), \label{eq:GProbEst}
\end{eqnarray}
where in the last step we observed that, given $\Phi\in \Xi(t)$, by the strong Markov property of $\xi$ applied successively to the stopping times
$\tau_i:=\inf\{j\geq \lfloor t/4\rfloor :\sum_{k=\lfloor t/4 \rfloor}^{j} 1_{A(t,\Phi,k, 0)}=i\}$, we can couple $\xi$ with a sequence
of i.i.d.\ Bernoulli random variables $(p_{i,\eps})_{i\in\N}$ with
$$
\P(p_{1,\epsilon}=1) = \P^\xi \big( Y_1^0(s) = 0 \, \forall s \in [0,\epsilon/\rho] \, \big \vert \, \xi(0,0) \geq 1 \big),
$$
such that $1_{A(t,\Phi,\tau_i, \eps)}\geq p_{i,\eps}$ a.s.\ for all $i\in\N$, and hence $\sum_{k \in \{\lfloor t/4 \rfloor,\cdots, t-1\}} 1_{A(t,\Phi,k, \epsilon)}\geq \sum_{i=1}^{t/2} p_{i,\eps}$ on the event $\sum_{k \in \{\lfloor t/4 \rfloor,\cdots, t-1\}} 1_{A(t,\Phi,k, 0)} \geq \frac{t}{2}$. Here $p_{i,\eps}$ corresponds to the event that given $A(t,\Phi, \tau_i, 0)$, a chosen $Y$-particle, which is close to $\Phi$ at time $\tau_i$, does not jump on the time interval $[\tau_i, \tau_i+\epsilon/\rho]$.

By \cite[Prop.~8]{KS05}, the first term in \eqref{eq:GProbEst} is bounded from above
by $1/t^{2}$ for $t$ large enough. For the second term we have $\big\vert \Xi (t)|_{\{\lfloor t/4 \rfloor, \cdots, t\}} \big\vert \leq e^{Ct}$
for some $C > 0$ and all $t,$
%
while large deviations yield that we can find $\epsilon > 0$ such that
$$
\P \Big( \sum_{k=1}^{t/2} p_{k,\epsilon} \leq \epsilon t \Big) \leq e^{-2Ct}
$$
for $t$ large enough. From now on we fix such an $\epsilon$.
Borel-Cantelli then yields that $\P^\xi$-a.s., $G(t)$ holds for all $t\in\N$ large enough.

Next observe that by the strong Markov property of $X$, we can construct a coupling such that
on the event $\sum_{k \in \{\lfloor t/4 \rfloor,\cdots, t-1\}} 1_{A(t,X,k, \epsilon)} \geq \epsilon t$, the random variable
$\int_0^t \xi(s, X(s)) \,{\rm d}s$ almost surely dominates the sum of i.i.d.\ random variables $(q_{i,\eps})_{1\leq i\leq \eps t}$
with
$$
\begin{aligned}
&\P(q_{1,\eps}=\eps/(2\rho)) &=&\ \alpha:= \inf_{y,z\in\mathcal{Q}_r(0)} \P^X_y(X(s)=z \, \forall s \in [\epsilon/(2 \rho),\epsilon/ \rho)]) >0, \\
&\P(q_{1,\eps}=0) &=&\ 1-\alpha;
\end{aligned}
$$
$q_{i,\eps}$ corresponds to the event that given $\tau_i:=\inf\{j\geq \lfloor t/4\rfloor :\sum_{k=\lfloor t/4 \rfloor}^{j} 1_{A(t,X,k, \eps)}=i\}$,
$X$ finds a $Y$-particle in the $\xi$ field which guarantees the event $A(t, X, \tau_i, \eps)$, and then occupies the same position as that $Y$-particle
on the time interval $[\tau_i+\eps/(2\rho), \tau_i+\eps/\rho]$. Since $\P^\xi$-a.s., $G(t)$ holds for all $t\in\N$ large enough, for such $t$, we have
\beq \label{eq:FKEst}
\E_0^X \Big[ \exp \Big\{ -\gamma \int_0^t \xi(s,X(s)) \, {\rm d}s \Big\}, \Xi(t) \Big]
\leq \E\Big[e^{-\gamma \sum_{i=1}^{\eps t} q_{i,\eps}}\Big] = \big( \alpha e^{-\gamma \epsilon/(2 \rho)} + 1-\alpha)^{\epsilon t}.
\eeq
Thus, with \eqref{eq:manyJumpsEst} and \eqref{eq:FKEst} we obtain that $\P^\xi$-a.s., for all $t\in\N$ large,
\begin{eqnarray*}
\E_0^X \Big[ \exp \Big\{ -\gamma \int_0^t \xi(s,X(s)) \, {\rm d}s \Big\} \Big]
&\leq & \E_0^X \Big[ \exp \Big\{ -\gamma \int_0^t \xi(s,X(s)) \, {\rm d}s \Big\}, \Xi(t) \Big]
+ \P_0^X \big( X \in \Xi(t)^c \big)\\
&\leq&  e^{-\delta(t+ o(t))}
\end{eqnarray*}
for some $\delta > 0$. This establishes the desired result along integer $t$. Since $Z_{t,\xi}^\gamma$ is monotone
in $t$, we deduce that the result holds as stated.
\qed
\bigskip

\noindent
{\bf Acknowledgement} We thank Frank den Hollander for bringing \cite{R94} to our attention, Alain-Sol Sznitman for suggesting that we
prove a shape theorem for the quenched survival probability, and Vladas Sidoravicius for explaining to us \cite[Prop.~8]{KS05}, which we use to prove
the positivity of the quenched Lyapunov exponent. A.F.~Ram\'irez was partially supported by Fondo Nacional de Desarrollo Cient\'\i fico y Tecnol\'ogico grant 1100298.
J.~G\"artner, R.~Sun and partially A.F.~Ram\'irez were supported by the DFG Forschergruppe 718 Analysis and Stochastics in Complex Physical Systems.

\end{document}